\documentclass[10pt,twoside]{article}
\usepackage{amsmath,amssymb,amscd,latexsym}
\usepackage[all]{xy}
\newcommand{\C}{{\mathbb{C}}}
\newcommand{\F}{{\mathbb{F}}}
\newcommand{\N}{{\mathbb{N}}}

\newcommand{\R}{{\mathbb{R}}}

\newcommand{\Z}{{\mathbb{Z}}}
\newcommand{\be}{\mathbf{1}}
\newcommand{\cpr}{\mathrm{cpr }\,}

\newcommand{\diam}{\mathrm{diam}}
\newcommand{\dist}{\mathrm{dist}}
\newcommand{\her}{\mathrm{her}}
\newcommand{\bi}{{\bar{\mbox{\it \i}}}}

\newcommand{\id}{\mathrm{id}}
\newcommand{\ord}{\mathrm{ord }\,}
\newcommand{\rk}{\mathrm{rk}\,}
\newcommand{\rr}{\mathrm{rr }\,}

\newcommand{\sr}{\mathrm{sr}\,}
\newcommand{\supp}{\mathrm{supp }\,}
\newcommand{\Sd}{\mathrm{Sd }\,}
\newcommand{\tr}{\mathrm{tr}\,}

\newcommand{\Bh}{{\mathcal B}}
\newcommand{\Ch}{{\mathcal C}}
\newcommand{\Fh}{{\cal F}}
\newcommand{\Hh}{{\mathcal H}}
\newcommand{\Kh}{{\mathcal K}}
\newcommand{\Oh}{{\mathcal O}}

\newcommand{\Th}{{\mathcal T}}
\newcommand{\Uh}{{\mathcal U}}
\newcommand{\Vh}{{\mathcal V}}

\newcommand{\tei}{\, | \,}

\newcommand{\verk}{\mbox{\scriptsize $\,\circ\,$}}
\newcommand{\halb}{\frac{1}{2}}
\newcommand{\te}{\textstyle}

\newcounter{number}[subsection]
\newcounter{altnumber}[section]
\newenvironment{nummer}{\refstepcounter{number}{\noindent\bf\arabic{section}.\arabic{subsection}.\arabic{number}}}{}
\newenvironment{altnummer}{\refstepcounter{altnumber}{\noindent\bf\arabic{section}.\arabic{altnumber}}}{}
\newcommand{\bn}{\noindent\begin{nummer} \rm}
\newcommand{\en}{\end{nummer}}

\newcommand{\altbn}{\noindent \begin{altnummer} \rm}
\newcommand{\alten}{\end{altnummer}}

\newenvironment{theorem}{\noindent {\bf Theorem:} \it}{}

\newenvironment{lemma}{\noindent {\bf Lemma:} \it}{}
\newenvironment{nlemma}{\noindent {\bf Lemma:} \it}{}
\newenvironment{prop}{\noindent {\bf Proposition:} \it}{}

\newenvironment{defn}{\noindent {\bf Definition:} \it}{}

\newenvironment{cor}{\noindent {\bf Corollary:} \it}{}
\newenvironment{ncor}{\noindent {\bf Corollary:} \it}{}
\newenvironment{conj}{\noindent {\bf Conjecture:} \it}{}
\newenvironment{remark}{\noindent {\bf Remark:}}{}
\newenvironment{nremark}{\noindent {\bf Remark:}}{}
\newenvironment{remarks}{\noindent {\bf Remarks:}}{}

\newenvironment{examples}{\noindent {\bf Examples:}}{}

\newenvironment{proof}{\noindent {\bf Proof:}}{\mbox{}\hfill$\Box$}

\begin{document}

\title{{\Large \sc Covering Dimension for\\ Nuclear $C^*$-Algebras }\vspace{3ex}}
\author{{\large Wilhelm Winter }\\ Mathematisches Institut\\ Universit\"at M\"unster\\ 
Einsteinstr. 62\\ 48149 M\"unster\\ Germany\\ e-mail: wwinter@uni-muenster.de }
\date{May 2001}
\maketitle
\setcounter{section}{-1}
\begin{abstract}

We introduce the completely positive rank, a notion of covering dimension for nuclear $C^*$-algebras and analyze some of its properties.

The completely positive rank behaves nicely with respect to direct sums, quotients, ideals and inductive limits. For abelian $C^*$-algebras it coincides with covering dimension of the spectrum and there are similar results for continuous trace algebras. 

As it turns out, a $C^*$-algebra is zero-dimensional precisely if it is $AF$. We consider various examples, particularly of one-dimensional $C^*$-algebras, like the irrational rotation algebras, the Bunce-Deddens algebras or Blackadar's simple unital projectionless $C^*$-algebra.

Finally, we compare the completely positive rank to other concepts of noncommutative covering dimension, such as stable or real rank.

\end{abstract}

\section[{\sc Introduction}]{\large \sc Introduction}

The theory of $C^*$-algebras is often regarded as noncommutative topology, mainly because any abelian $C^*$-algebra is completely determined by its spectrum, a locally compact space. So $C^*$-algebras might be thought of as noncommutative topological spaces and this point of view has been very fruitful, especially since the introduction of methods from algebraic topology (such as $Ext$ or $K$-theory) to $C^*$-algebra theory.

It is thus natural to look for invariants of topological spaces which can be defined analogously for $C^*$-algebras. A good candidate for such an invariant certainly is topological covering dimension and there have been several successful approaches in this direction, especially the stable and the real rank, introduced by Rieffel and by Brown and Pedersen, respectively. It is the aim of these notes to introduce another version of covering dimension, this time for nuclear $C^*$-algebras.

It is well-known that a $C^*$-algebra $A$ is nuclear if and only if it has the completely positive approximation property, i.e.\ the identity map on $A$ can be approximated by completely positive contractions of the form $ A \stackrel{\psi}{\longrightarrow} F \stackrel{\varphi}{\longrightarrow} A$ with finite-dimensional $C^*$-algebras $F$. By a system of c.p.\ approximations we then mean a net $((F_\lambda, \psi_\lambda, \varphi_\lambda))_\Lambda$ with $ A \stackrel{\psi_\lambda}{\longrightarrow} F_\lambda \stackrel{\varphi_\lambda}{\longrightarrow} A$ as above such that $\varphi_\lambda \verk \psi_\lambda \rightarrow \id_A$ pointwise.

An important step towards a better understanding of such systems of c.p.\ approximations certainly is the notion of generalized inductive limits and (strong) $NF$ algebras introduced by Blackadar and Kirchberg. The present paper might also be viewed as a step in that direction, for we use special structural properties of c.p.\ approximations to obtain a topological invariant.

The triples $(F_k, \psi_k, \varphi_k)$ may be thought of as analogues of open coverings of spaces and, imposing a certain condition on the $\varphi_k$, can be used to define the completely positive rank of $A$, $\cpr A$, our notion of noncommutative covering dimension.

As one would expect from a generalized dimension theory, for commutative $C^*$-algebras ${\cal{C}}_0(X)$ we have $\cpr({\cal{C}}_0(X)) = \dim X$, where $\dim X$ denotes ordinary covering dimension of the spectrum. The completely positive rank behaves nicely with respect to direct sums, quotients, ideals  and inductive limits, i.e.\ $\cpr(A \oplus B) = \max(\cpr A, \cpr B)$, $\cpr(A/J) \le \cpr A$, $\cpr J \le \cpr A$ and $\cpr(\lim_\rightarrow A_k) \le \underline{\lim} \, \cpr A_k$. 

$AF$ algebras, which are inductive limits of finite-dimensional (in the vector space sense) $C^*$-algebras, are easily seen to have completely positive rank zero, but it turns out that conversely $\cpr A = 0$ implies that $A$ is $AF$.

It is not hard to see that $AT$ algebras have completely positive rank smaller than one. In particular Bunce-Deddens algebras and irrational rotation algebras are on-dimensional, just as Blackadar's simple unital projectionless $C^*$-algebra.

As we pointed out before, the completely positive rank of a commutative $C^*$-algebra is equal to the covering dimension of its spectrum. One might ask if this also is the case in more general situations, e.g. for continuous trace algebras. It turns out that in fact $\cpr A \le \dim \hat{A}$ for any continuous trace algebra $A$. If the dimensions of the irreducible representations of $A$ are (at least locally) bounded, then for the strong completely positive rank we even have equality.

We also compare the completely positive rank to other concepts of noncommutative covering dimension such as stable, real, analytic and tracial rank. Although we do not have general results, our examples show that the completely positive rank does not coincide with any of these other concepts; it might well be that it dominates the real and the stable rank. 

The paper is organized as follows. First we give a brief survey of completely positive maps and nuclear $C^*$-algebras. In Section 1.3 we have collected a number of technical results on positive elements and projections which will be needed later. 

Section 2 recalls the notion of topological covering dimension and some several well-known results. We slightly modify the notion of the order of an open covering and give another characterization of covering dimension, which motivates the introduction of the completely positive rank in Section 3.

There we note some basic features of the theory and show that, for commutative $C^*$-algebras, it coincides with covering dimension of the spectrum.

Section 4 is concerned with the zero-dimensional case; 4.1 contains an explicit description of maps of strict order zero and in 4.2 we prove that a $C^*$-algebra has completely positive rank zero if and only if it is $AF$.

Next we recall some facts on continuous trace algebras and show that $\cpr A \le \dim \hat{A}$ for such $C^*$-algebras.

The completely positive rank is compared to other concepts of non-commutative covering dimension in Section 6.

Finally, in Section 8 we note some open questions related to the completely positive rank and very roughly outline possible variations of the concept.

I would like to thank my advisor, Prof. J. Cuntz, for drawing my attention to this subject and for many useful comments. I also thank Prof. B. Blackadar, Prof. S. Echterhoff and Prof. N. C. Phillips for many hints and discussions. Furthermore, a big ``Dankesch\"on'' to H. Eissing for providing excellent working conditions. I am gratefully indebted to G. Weckermann, who converted an almost unreadable manuscript into this paper (it was hard to convince her that not all nuclear $C^*$-algebras are unclear). Finally, it is a pleasure to thank my friend and colleague B. Neub\"user for his invaluable help in dealing with unfriendly computers and, of course, for lots of fruitful discussions.
\section[{\sc $C^*$-algebra preliminaries}]{\large \sc $C^*$-algebra preliminaries}

\subsection[{\rm Completely positive maps}]{\large \sc Completely positive maps}

For the convenience of the reader we give a brief survey of completely positive maps. A good general reference is \cite{Pa}.

\bn
A linear map $\varphi : A \to B$ between $C^*$-algebras is called positive, if $\varphi (A_+) \subset B_+$; it is then automatically $*$-preserving. $\varphi$ is $n$-positive, if the induced map $\varphi^{(n)} : M_n (A) \to M_n (B)$ is positive, and completely positive (c.p.), if it is $n$-positive for all $n \in \N$. If $A$ and $B$ are unital and $\varphi (\be_A) = \be_B$, we say $\varphi$ is unital completely positive (u.c.p.).

It makes sense (and will be useful) to allow also completely positive maps between operator systems, i.e.\ unital self-adjoint linear subspaces of unital $C^*$-algebras (or even between operator spaces, i.e.\ linear subspaces of $C^*$-algebras).

Clearly, all $*$-homomorphisms are completely positive. If $\varphi \, : \, A \rightarrow  B$ is c.p.\ and $h \in B$, then the map $\varphi_h (\,.\,) := h^* \, \varphi (\,.\,) \, h$ is c.p.\ as well. One has $\| \varphi \| = \| \varphi^{(n)} \| = \lim_{\lambda} \| \varphi (u_{\lambda}) \|$, where $(u_{\lambda})_{\Lambda}$ is some positive approximate unit for $A$ with $\|u_{\lambda} \| \le 1 \; \forall \lambda \in \Lambda$. We will deal almost exclusively with contractive (i.e.\ norm-decreasing) c.p.\ maps.

\en

\bn{\label{cp-unitization}}
Many results on completely positive maps in the literature are stated for unital maps, but one can often give nonunital versions using the following facts:

\begin{nlemma} {\rm{(\cite{CE1}, Lemma 3.9)}}
Let $A,B$ be $C^*$-algebras and $\varphi : A \to B$ a completely positive contraction. Then $\varphi$ extends (uniquely) to a unital completely positive map $\varphi^+ : A^+ \to B^+$.
\end{nlemma}
\en

\bn{\label{unital-compression}}
\begin{lemma} {\rm{(\cite{CE2}, Lemma 2.2)}}
Let $X$ be an operator system in some $C^*$-algebra, $N$ a von Neumann algebra and $\varphi : X \to N$ completely positive. Then there is a unital completely positive map $\sigma : X \to N$ such that $\varphi (a) = \varphi(\be)^{\halb} \sigma (a) \varphi (\be)^{\halb} \; \forall a \in X$.
\end{lemma}
\en

\bn{\label{Stinespring}}
The importance of completely positive maps largely comes from their characterization as compressions of $*$-homomorphisms given by Stinespring. We state Lance's version of this theorem (\cite{La1}, Theorem 4.1), which does not assume $A$ to be unital. The proof is a generalization of the  GNS--construction.

\begin{theorem}
Let $A$ be a $C^*$-algebra and $\varphi : A \to \Bh (\Hh)$ a completely positive map.\\
Then there are a Hilbert space $\Hh'$, a $*$-homomorphism $\pi : A \to \Bh (\Hh')$ and a bounded operator $V : \Hh \to \Hh'$ such that
\[
\varphi (a) = V^* \pi (a) V \; \forall a \in A.
\]
\end{theorem}

\begin{nremark}
  If $\varphi$ is unital, $V$ is an isometry and $\Hh$ may be regarded as a subspace of $\Hh'$. With this identification we have $\varphi (a) = p_{\Hh} \pi (a) p_{\Hh} $, so indeed $\varphi$ is a compression of a $*$-homomorphism.
\end{nremark}
\en

\bn
\begin{ncor}
{\rm{(cf.\ \cite{BK1}, Corollary 4.1.3)}}
Let $A$, $B$ be $C^*$-algebras and $\varphi : A \to B$ completely positive contractive.\\
For all $x \in A$ we have $\varphi (x^* x) \ge \varphi (x)^* \varphi (x)$. If $\varphi (x^* x) = \varphi (x)^* \varphi (x)$, then $\varphi (yx) = \varphi (y) \varphi (x)$ for all $y \in A$.
\end{ncor}
\en

\bn{\label{multiplicativity}}
\begin{prop}
Let $A$, $B$ be $C^*$-algebras, $\varphi : A \to B$ completely positive contractive, $x \in A$ with $\| x \| \le 1$ and such that $\| \varphi (x^* x) - \varphi (x)^* \varphi (x) \| \le \varepsilon$ for some $\varepsilon \ge 0$.\\
Then $\| \varphi (yx) - \varphi (y) \varphi (x) \| \le \sqrt{\varepsilon}$ for all $y \in A$ with $\| y \| \le 1$.
\end{prop}

\begin{proof}
  This is a variation of \cite{BK1}, Corollary 4.1.5. We may assume $A$ and $\varphi$ to be unital and $B = \Bh (\Hh)$. By Stinespring's theorem we have that $\varphi (a) = p \pi (a) p  \; \forall \, a \in A$, where $\pi : A \to \Bh (\Hh')$ is a $*$-representation and $p \in \Bh (\Hh')$ is the projection onto $\Hh \subset \Hh'$, i.e.\ $\Bh (\Hh) = p \Bh (\Hh')p$.\\
Then we get
\[
\| \pi (x) p - p \pi (x) p \|^2 = \| p \pi (x^* x) p - p \pi (x^*) p \pi (x) p \| = \| \varphi (x^* x) - \varphi (x)^* \varphi(x) \| \le \varepsilon \; ,
\]
hence $\| \pi (x) p - p \pi (x) p \| \le \sqrt{\varepsilon}$. Finally
\[
\| \varphi (yx) - \varphi (y) \varphi (x) \| = \| p \pi (y) \pi (x) p - p \pi (y) p \pi (x) p \| \le \sqrt{\varepsilon} \; .
\]
\end{proof}

\begin{remark}
Note that the condition of the preceding proposition is fulfilled, if $x \ge 0$ and $\| \varphi (x) - \varphi (x)^2 \| \le \varepsilon$, since $0 \le \varphi (x^2) - \varphi (x)^2 \le \varphi (x) - \varphi (x)^2$.
\end{remark}
\en

\bn{\label{Arveson}}
We will also need Arveson's extension theorem, which says that, for any Hilbert space $\Hh$, $\Bh (\Hh)$ is injective in the category of operator systems with completely positive maps.

\begin{theorem}
Let $Y \subset X$ be operator systems in some $C^*$-algebra, $\varphi : Y \to \Bh (\Hh)$ a completely positive map.\\
Then $\varphi$ extends to a completely positive map $\tilde{\varphi} : X \to \Bh (\Hh)$.
\end{theorem}
\en

\subsection[{\rm Nuclear $C^*$-algebras}]{\large \sc Nuclear $C^*$-algebras}
We recall the notion of nuclear $C^*$-algebras and some structure results. For more information and for attributions see \cite{La2} or \cite{Wa}.

\bn
In general there are many $C^*$-norms on the algebraic tensor product $A \odot B$ of two $C^*$-algebras but there is always a maximal and a minimal one. $A$ is called nuclear, if for any $B$ there is only one $C^*$-norm on $A \odot B$, i.e.\ $\| \cdot \|_{\max} = \| \cdot \|_{\min}$.

Many of our stock-in-trade $C^*$-algebras are nuclear:
\begin{itemize}
\item finite-dimensional $C^*$-algebras
\item abelian $C^*$-algebras
\item inductive limits of nuclear $C^*$-algebras, in particular $AF$ algebras
\item continuous trace $C^*$-algebras
\item extensions of nuclear by nuclear $C^*$-algebras
\item ideals and quotients of nuclear $C^*$-algebras
\item tensor products of nuclear $C^*$-algebras
\end{itemize}
are all nuclear.
\en

\bn{\label{cpap}}
\begin{defn}
We say that a $C^*$-algebra $A$ has the completely positive approximation property, if the following holds:\\
For $a_1 , \ldots , a_k \in A$ and $\varepsilon > 0$ there are a finite-dimensional $C^*$-algebra $F$ and completely positive contractions $\psi : A \to F$ and $\varphi : F \to A$ such that 
\[
\| \varphi \psi (a_i) - a_i \| < \varepsilon \; \forall i \; .
\]
We call the triple $(F , \psi , \varphi)$ a c.p.\ approximation within $\varepsilon$ for $\{ a_1 , \ldots , a_k \}$. A system of c.p.\ approximations for $A$ is a net $((F_{\lambda} , \psi_{\lambda} , \varphi_{\lambda}))_{\Lambda}$ of c.p.\ approximations such that $\varphi_{\lambda} \psi_{\lambda} \to \id_A$ pointwise. 
\end{defn}

We will almost exclusively deal with separable $C^*$-algebras, thus it will always suffice to consider systems with index set $\Lambda = \N$. 
\en

\bn{\label{rcpap}}
\begin{remarks} \rm
(i) This definition will be crucial for our notion of noncommutative covering dimension: In Definition \ref{cpr} we will replace open coverings of order $n$ (cf.\ Definition \ref{covering-dimension}) by c.p.\ approximations $(F, \psi , \varphi)$ together with an extra condition on $\varphi$.\\
(ii) Instead of studying a nuclear $C^*$-algebra it is as good to analyze the properties of its systems of c.p.\ approximations. One approach in this direction is the notion of generalized inductive limits (and that of (strong) $NF$ algebras) introduced by Blackadar and Kirchberg in \cite{BK1}. \\
The present paper should also be viewed as an attempt to study the fine structure of nuclear $C^*$-algebras in terms of systems of c.p.\ approximations.
\end{remarks}
\en

\bn{\label{tcpap}}
We then have the following characterization of nuclearity (cf.\ \cite{Wa}, Propositions 2.1 and 2.2):

\begin{theorem}
A $C^*$-algebra $A$ is nuclear if and only if it has the completely positive approximation property.
\end{theorem}
\en

\bn{\label{lifting}}
A $C^*$-algebra $A$ is said to have the lifting property, if for any $C^*$-algebra $B$ containing an ideal $J$ every completely positive contraction $\varphi : A \to B / J$ has a c.p.c.\ lifting $\psi : A \to B$, i.e.\ one has a commuting diagram
\[
\xymatrix{
 & B \ar[d]\\
A \ar@{-->}[ur]^{\exists \, \psi} \ar[r]^{\varphi} & B/J \; .
}
\]

\begin{theorem} {\rm{(Choi-Effros Lifting Theorem, \cite{CE1}, Theorem 3.10)}} 
Every separable nuclear $C^*$-algebra has the lifting property.
\end{theorem}
\en

\subsection[{\rm Some technicalities}]{\large \sc Some technicalities}
\bn{\label{functions}}
We first introduce some notation which is used frequently throughout the paper. For positive numbers $\alpha$ and $\varepsilon$ define positive functions in $\Ch (\R)$ as follows:
\[
f_{\alpha,\varepsilon} (t) :=  \left\{ 
    \begin{array}{cl}
0 & \mbox{for} \; t \le \alpha \\
t & \mbox{for} \; \alpha + \varepsilon \le t \\
\mbox{linear} &  \mbox{elsewhere}
\end{array} \right.
\]
and
\[
g_{\alpha , \varepsilon} (t) :=  \left\{ 
    \begin{array}{cl}
0 & \mbox{for} \; t \le \alpha \\
1 & \mbox{for} \; \alpha + \varepsilon \le t \\
\mbox{linear} & \mbox{elsewhere \,;}
\end{array} \right.
\]
we write $g_{\alpha}$ for the characteristic function of $[\alpha , \infty]$. If $h$ is a positive element in some $C^*$-algebra $A$, then $\| f_{\varepsilon , \varepsilon} (h) - h \| \le \varepsilon$. If the spectrum $\sigma (h)$ is disconnected, say $\sigma (h) \subset [0, \varepsilon] \cup [1 - \varepsilon , 1]$ for $\varepsilon < \halb$, then $g_{\halb} (h)$ is a projection. For the moment set
\[
f^- (t) := 
\left\{ 
\begin{array}{cl}
0 & \mbox{for} \; t \le \varepsilon \\
t^{-1} & \mbox{for} \; 1 - \varepsilon \le t \\
\mbox{linear} & \mbox{elsewhere \,,}
\end{array} \right.
\]
then $f^- (h) \, h = h \, f^- (h) = g_{\halb} (h)$. In this case we also write $h^{-1}$ for $f^- (h)$ to indicate that $f^- (h)$ is the inverse of $h$, taken in $g_{\halb} (h) A g_{\halb} (h)$.
\en 

\bn
Recall that elements $x,y$ in a $C^*$-algebra $A$ are said to be orthogonal, $x \perp y$, if $xy = yx = x^* y = x y^* = 0$. Subsets $X,Y \subset A$ are orthogonal, $X \perp Y$, if $x \perp y$ for all $x \in X$ and $y \in Y$.\\
If $x,y \in A_{sa}$, then $xy = 0 \Leftrightarrow yx = 0$. If $x,y \in A_+$ and $0 \le x' \le x , \  0 \le y' \le y$, then $x \perp y \Rightarrow x' \perp y'$. 
\en

\bn
Let $B$ be a $C^*$-subalgebra of the $C^*$-algebra $A$. $B$ is said to be a hereditary subalgebra, $B \subset_{\her} A$, if 
\[
0 \le x \le y , \ y \in B \Longrightarrow x \in B.
\]
Any ideal in $A$ automatically is hereditary. If $B$ has a strictly positive element $h$ (or, equivalently, if $B$ is $\sigma$-unital), then $B = \overline{hBh} = \overline{hAh}$. On the other hand, if $h \in A_+$ is any positive element, then $\overline{hAh} \subset_{\her} A$ and $h$ is strictly positive for $\overline{hAh}$. We write $A_h$ for $\overline{hAh}$; note that, if $h$ happens to be a projection, then $\overline{hAh} = hAh$.\\
If $a,b \in A_+$, then $a \perp b$ if and only if $A_a \perp A_b$.
\en

\bn{\label{almost-units}}
\begin{prop}
Let $A$ be a $C^*$-algebra, $d \ge h \in A_+ , \, \| d \| \le 1 , \, x \in A$ with $\| x \| \le 1$ and $\varepsilon \ge 0$. If $\| (\be -h) x \| \le \varepsilon$, then $\| (\be -d) x \| \le \sqrt{\varepsilon}$.
\end{prop}

\begin{proof}
  $\| (\be-d) x \|^2 = \| x^* (\be -d)^2 x \| \le \| x^* (\be-d) x \| \le \|x^* (\be-h) x \| \le \varepsilon$.
\end{proof}
\en

\bn{\label{almost-projections}}
\begin{prop}
Let $A$ be a $C^*$-algebra and $h \in A_+ , \, \| h \| \le 1$, such that $\| h - h^2 \| < \varepsilon < \frac{1}{4}$.\\
Then there is a projection $p \in C^* (h) \subset A$ with $\| p - h \| < 2 \varepsilon$. Furthermore, $php$ is invertible in $p C^* (h) p$ and for $c := (php)^{-\halb}$ one has $\| p-c \| < 4\varepsilon$.
\end{prop}

\begin{proof}
  This is an argument used frequently in $K$-theory, cf.\ e.g.\ \cite{WO}, Lemma 5.1.6. If $t \in \sigma (h)$, then $0 \le t - t^2 < \varepsilon < \frac{1}{4}$, so $\sigma (h) \subset [ 0 , \halb - \delta] \cup [\halb + \delta , 1]$, where $\delta := \halb \sqrt{1 - 4 \varepsilon} > 0$.\\
But then $g_{\halb}$ is continuous on $\sigma (h)$ and $p := g_{\halb} (h)$ is a projection satisfying 
\[
\| p - h \| = \| g_{\halb} (h) - \id (h) \| \le \halb - \delta < 2 \varepsilon
\]
(see \ref{functions} for the definition of $g_{\alpha}$). Finally, $(php)^{-\halb} = f (h)$, where
\[
f (t) := \left\{ 
  \begin{array}{ccl}
0 & \mbox{for} & t \le \halb - \delta \\
t^{-\halb} & \mbox{for} & t \ge \halb + \delta \; ,
  \end{array} \right.
\]
and one has
\[
\| p - c \| = \| g_{\halb} (h) - f (h) \| \le \frac{1}{\halb + \delta} - 1 \le 2 \cdot \left( \halb - \delta \right) < 4 \varepsilon \; .
\]
\end{proof}
\en

\bn{\label{orthogonal-projections}}
\begin{cor}
Let $A$ be a $C^*$-algebra, $0 \le \delta < \frac{1}{24}$ and $p,q \in A$ projections satisfying $\| pq \| \le \delta$.\\
Then there is a projection $\tilde{p} \in A$ with $\tilde{p} \perp q$ and $\| \tilde{p} - p \| \le 14\delta$.
\end{cor}

\begin{proof}
  Set $h := (\be-q) \, p \, (\be-q)$, then
\[
\| h-p \| = \| qpq - qp - pq \| \le \| qp(\be-q) \| + \|pq\| \le 2\delta \; .
\]
Furthermore we have
\[
\| h^2 - h \| \le \| h^2 - ph \| + \| ph - p \| + \| p-h \| \le 6 \delta \; ,
\]
and by Proposition \ref{almost-projections} there exists a projection $\tilde{p} \in C^* (h) \subset A$ satisfying $\| \tilde{p} - h \| \le 12 \delta$, thus $\| \tilde{p} - p \| \le 14 \delta$. Obviously $\tilde{p} \perp q$.
\end{proof}
\en

\bn{\label{close-projections}}
\begin{prop}
Let $p,q$ be projections in a $C^*$-algebra $A$ with $\| p-q \| < \eta \le \frac{1}{4}$. Then there is a partial isometry $s \in A$ such that $s^*s = p , ss^* = q$ and $\| s-p \| < 4\eta$.
\end{prop}

\begin{proof}
%
%
This basically is \cite{Cu3}, Lemma 2.2. If necessary, we adjoin a unit to $A$.\\
Define symmetries $x := 2 p - \be$ and $y := 2 q - \be$, then there is $h \in A_{sa}$ with $\sigma (h) \subset (-\pi , \pi)$ such that $xy = e^{ih}$. Also, one readily checks that $\| \be - xy \| < 4 \eta$.\\
Set $u:= e^{-\frac{ih}{2}}$, then it is not hard to see that $\| u - \be \| < 4 \eta$ as well. As in \cite{Cu3}, Lemma 2.2, one checks that $u x u^* = y$, thus $u p u^* = q$ and $u^* q u = p$.\\
Set $s := u p \in A$, then $s^* s = p$, $s s^* = q$ and 
\[
\|s - p\| \le \|u - \be \| < 4 \eta \, .
\]
\end{proof}
\en

We close with some simple technical observations which will be useful in Section 3. By $U(M_r)$ we denote the unitary group of $M_r$.

\bn{\label{unitary-neighborhoods}}
\begin{lemma}
(i) Let $e_1 , \ldots , e_r \in M_r$ be minimal projections such that the corresponding unit vectors $\xi_i \in \C^r$ are linearly independent.
Then $e_1 + \ldots + e_r$ is (positive and) invertible in $M_r$, i.e.\ $\be_{M_r} \le \lambda \cdot (e_1 + \ldots + e_r)$ for some $\lambda > 0$.\\
(ii) If $p \in M_r$ is a minimal projection and $\emptyset \neq \Uh \subset U (M_r)$ an open subset of $U (M_r)$, then there are unitaries $u_1 , \ldots , u_r \in \Uh$ such that $u^*_1 p u_1 + \ldots + u^*_r p u_r$ is invertible in $M_r$.
\end{lemma}

\begin{proof}
  (i) If $e_1 + \ldots + e_r$ is not invertible, then there is some $\eta \in \C^r$ such that (in $\langle \mbox{bra}|$-$|\mbox{ket} \rangle$ notation),
\[
0 = \langle \eta | e_1 + \ldots + e_r | \eta \rangle = \langle \eta |e_1| \eta \rangle + \ldots + \langle \eta |e_r| \eta \rangle.
\]
But then $\langle \eta |e_i| \eta\rangle = 0 \; \forall \, i$; in particular $\eta$ is orthogonal to all the $\xi_i$. We would thus have $r+1$ linearly independent vectors in $\C^r$, a contradiction.\\
(ii) follows from (i), for, since $\Uh$ is open, one can choose unitaries $u_1 , \ldots , u_r$ in $\Uh$ such that $u_1 \xi , \ldots , u_r \xi$ are linearly independent (where $\xi$ denotes a unit vector corresponding to $p$).
\end{proof}
\en

\bn{\label{dominated-orthogonality}}
\begin{remark}
We shall use Lemma \ref{unitary-neighborhoods} frequently in the following context: 

Let $A, \, F$ be $C^*$-algebras, $F$ finite-dimensional and $\varphi : F \to A$ be c.p.c.\ Suppose $M_r$ sits in $F$ as a hereditary subalgebra and $\emptyset \neq \Uh \subset U(M_r)$ is open.\\
Consider minimal projections $e, \, \bar{e} \in F$ such that $e \in M_r$ and $\bar{e} \perp M_r$. Then the lemma implies that there are unitaries $u_1, \ldots , u_r \in \Uh$ such that $h := u_1^* e u_1 + \ldots + u_r^* e u_r$ is invertible in $M_r$, hence $\be_{M_r} \le \lambda \cdot h$ for some $\lambda > 0$.\\
As a consequence, if $\varphi (u^* e u) \perp \varphi (\bar{e}) \; \forall u \in \Uh$, we have $\varphi (\be_{M_r}) \perp \varphi (\bar{e})$.
\end{remark}
\en


\section[{\sc Topological covering dimension}]{\large \sc Topological covering dimension}

In this section we recall the notion of topological covering dimension, describe some of the most important properties and give an equivalent characterization. See \cite{Pe} and \cite{En} for further information. 

\altbn{\label{covering-dimension}}
Recall the definition of the covering dimension of a topological space:

\begin{defn}
{\rm{(\cite{Pe}, Definition 3.1.1)}}
Let $X$ be a topological space.\\
a) The order of a family $(U_{\lambda})_{\Lambda}$ of subsets of $X$ does not exceed $n$, if for any $n +2 $ distinct indices $\lambda_0 , \ldots , \lambda_{n+1} \in \Lambda$ we have $\bigcap^{n+1}_{i=0} U_{\lambda_i} = \emptyset$.\\
b) The covering dimension of $X$ does not exceed $n$, $\dim X \le n$, if every finite open covering of $X$ has an open refinement of order not exceeding $n$. We say $\dim X = n$, if $n$ is the least integer such that $\dim X \le n$.
\end{defn}
\alten

\altbn{\label{dim-properties}}
One of the original motivations for dimension theory was to find a distinguishing invariant for the Euclidean spaces $\R^n$; in fact one has $\dim \R^n = n$.

Let $K \subset X$ be a closed subset, then $\dim K \le \dim X$. If $X$ is, say, compact and metrizable, we even have $\dim X = \max \{ \dim X , \, \dim (X \setminus K) \}$. For nonempty compact spaces we have $\dim (X \times Y) \le \dim X + \dim Y$.
\alten

\altbn{\label{countable-sums}}
We shall also need the following result, often referred to as the countable sum theorem:

\begin{theorem} {\rm (\cite{Pe}, Theorem 3.2.5)}
Let $X$ be a normal space such that $X = \bigcup_{i \in \N} A_i$ for closed subsets $A_i$ with $\dim A_i \le n$.\\
Then $\dim X \le n$.
\end{theorem}
\alten

\altbn{\label{limits-dim}}
\begin{prop}
{\rm{(\cite{Pe}, Corollary 8.1.7)}}
Let $X = \lim_{\leftarrow} X_{\alpha}$ for an inverse system $(X_{\alpha} , \pi_{\alpha,\beta})_{\alpha,\beta \in \Omega}$ of compact spaces $X_{\alpha}$ with $\dim X_{\alpha} \le n$ for each $\alpha \in \Omega$.\\
Then $X$ is compact and $\dim X \le n$.
\end{prop}
\alten

\altbn{\label{dim=0}}
There is a nice characterization of zero-dimensional locally compact spaces (cf. \cite{En}, Theorem 1.4.5):

\begin{prop}
For a locally compact space $X$, the following are equivalent:\\
(i) $\dim X = 0$,\\
(ii) $X$ is totally disconnected, i.e.\ it has no connected subspace consisting of more than one point.
\end{prop}
\alten

\altbn
{\bf Example:} A typical example of a zero-dimensional space is the Cantor set $C$, which may be obtained by the classical middle thirds construction or as the cartesian product $\{ 0,1 \}^{\N}$. It is typical in the sense that it is a universal space for all zero-dimensional second countable locally compact spaces, i.e.\ every such space can be embedded into $C$ (cf. \cite{En}, Theorem 1.3.15).
\alten

\altbn{\label{d-strict-order}}
Our goal is to define an analogue of covering dimension for nuclear $C^*$-algebras, using c.p.\ approximations. We could do this along the lines of Definition \ref{covering-dimension}, but the resulting condition on the approximations would be somewhat unnatural. We therefore introduce another characterization of covering dimension, which is certainly well-known to topologists, although we could not find it in the literature.

\begin{defn}
The strict order of a family $(U_{\lambda})_{\Lambda}$ of subsets of a topological space $X$ does not exceed $n$, if for any $n+2$ distinct indices $\lambda_0 , \ldots , \lambda_{n+1} \in \Lambda$ there are $i,j \in \{ 0 , \ldots , n+1 \}$ such that $U_{\lambda_i} \cap U_{\lambda_j} = \emptyset$.
\end{defn}
\alten

\altbn{\label{p-strict-order}}
\begin{prop}
Let $X$ be a normal space. T.f.a.e.:\\
i) $\dim X \le n$\\
ii) every finite open covering of $X$ has a (finite) open refinement of strict order not exceeding $n$.
\end{prop}

\begin{proof}
  b) $\Rightarrow$ a) is immediate. We show a) $\Rightarrow$ b):\\
Take a finite open covering of $X$. Because $\dim X \le n$, there is an open refinement $(U'_{\lambda})_{\Lambda}$ of order not exceeding $n$. By a standard argument (cf. \cite{Pe}, Proposition 3.1.2) we may assume $\Lambda$ to be finite and set $k := |\Lambda|$.\\
Take a partition of unity $(h_{\lambda})_{\Lambda}$ subordinate to $(U'_{\lambda})_{\Lambda}$ and set 
\[
U_{\lambda} := h^{-1}_{\lambda} ((0,1]) \, ,
\]
then $(U_{\lambda})_{\Lambda}$ also is a refinement of order not exceeding $n$. It suffices to find an open refinement $(V_{\gamma})_{\Gamma}$ of $(U_{\lambda})_{\Lambda}$ of strict order not exceeding $n$.\\
View $(h_{\lambda})_{\Lambda}$ as map $h : X \to K \subset \Delta_{k-1} \subset \R^k$, where $\Delta_{k-1}$ is the standard simplex in $\R^k$ and $K \subset \Delta_{k-1}$ is the minimal subcomplex of $\Delta_{k-1}$ containing $h (X)$ (s. \cite{ES}, Ch. 2 or \cite{Pe}, 2.6 for an introduction to simplicial complexes).\\
Let $u_{\lambda}$ be the vertices of $\Delta_{k-1}$ and $A_{\lambda}$ be the open stars around $u_{\lambda} , \lambda \in \Lambda$. Then $U_{\lambda} = h^{-1} (A_{\lambda})$;  as $(U_{\lambda})_{\Lambda}$ has order less than or equal to $n$, we have $\dim K \le n$ (note that $\dim K = \dim \tilde{K}$, where $\dim \tilde{K}$ is the (combinatorial) dimension of the (abstract) simplicial complex $\tilde{K}$, the geometric realization of which is $|\tilde{K}| \approx K$).

Let now $\Sd K$ be the barycentric subdivision of $K$ with vertices $v_{\gamma}$ and open stars $B_{\gamma} , \gamma \in \Gamma$ (cf. \cite{ES} 2.6).\\
Set $V_{\gamma} := h^{-1} (B_{\gamma})$, then $(V_{\gamma})_{\Gamma}$ is a refinement of $(U_{\lambda})_{\Lambda} ; \:(V_{\gamma})_{\Gamma}$ is of order less than or equal to $n$, as $\dim (\Sd K) = \dim K \le n$. Note that $(V_{\gamma})_{\Gamma}$ is of strict order not exceeding $n$ if the following condition on $\Sd K$ holds:
\begin{itemize}
\item[($\ast$)] for any distinct vertices $v_{\gamma_0} , \ldots , v_{\gamma_{n+1}}$ of $\Sd K$ there exist indices $i,j \in \{ 0 , \ldots , n+1 \}$ such that $v_{\gamma_i}$ and $v_{\gamma_j}$ do not sit in one and the same face of $\Sd K$, i.e.\ $v_{\gamma_i}$ and $v_{\gamma_j}$ are not connected by an edge of $\Sd K$.
\end{itemize}
The vertices $v_{\gamma}$ of $\Sd K$ are exactly the barycenters of simplexes $s_{\gamma}$ of $K$. But then $v_{\gamma_0} , \ldots , v_{\gamma_q}$ span a simplex in $\Sd K$, if $s_{\gamma_i}$ is a face of $s_{\gamma_{i+1}}$ in $K$ for $i = 0 , \ldots , q-1$, and every simplex in $\Sd K$ is of this form (cf. \cite{ES} 2.6).

Let now $v_{\gamma_0} , \ldots , v_{\gamma_q}$ be distinct vertices of $\Sd K$, such that $v_{\gamma_i}$ is connected to $v_{\gamma_j}$ by an edge in $\Sd K$ if $i \neq j$. But then $\dim s_{\gamma_i} \neq \dim s_{\gamma_j}$, because $s_{\gamma_i}$ is a proper face of $s_{\gamma_j}$ or vice versa. \\
As $0 \le \dim s_{\gamma} \le n \;\; \forall \gamma$, we have $q \le n$. Thus ($\ast$) holds and we are done.
\end{proof}
\alten

\section[{\sc Completely positive rank}]{\large \sc Completely positive rank}

In this section we define the completely positive rank and deduce some basic properties. In particular we show that, for commutative $C^*$-algebras, our theory coincides with covering dimension of the spectrum.

\altbn{\label{elementary}}
As we already pointed out, we regard a c.p.\ approximation $(F,\psi,\varphi)$ for a nuclear $C^*$-algebra $A$ as an analogue of an open covering. The strict order of such a c.p.\ approximation is then expressed as a condition on $\varphi$:

\begin{defn}
Let $A,F$ be $C^*$-algebras, $F$ finite-dimensional.\\
a) We say a set $\{ e_0 , \ldots , e_n \} \subset F$ is elementary, if the $e_i$ are mutually orthogonal minimal projections.\\
b) A completely positive map $\varphi : F \to A$ is of strict order not exceeding $n , \ord \varphi \le n$, if the following holds:\\
For every elementary set $\{ e_0 , \ldots , e_{n+1} \} \subset F$ there exist $i,j \in \{ 0 , \ldots , n+1 \}$ such that $\varphi(e_i) \perp \varphi(e_j)$. 
\end{defn}
\alten

\altbn{\label{strict-order}}
\begin{remarks}
Let $F$ be finite-dimensional, $A = \Ch_0 (X)$ for some locally compact space $X$ and $\varphi : F \to A$ a c.p.\ contraction.\\
(i) If  $F = \C^k$ for some $k$, then we have $\ord \varphi \le n$ if and only if the family of open sets 
\[
(\varphi (e_i)^{-1} (( 0, \infty)))_{i = 1 , \ldots , k}
\]
is of strict order $n$ in the sense of Definition \ref{d-strict-order}.\\
(ii) If $\ord \varphi \le n$, then for every set $\{p_0, \ldots , p_{n+1} \} \subset F$ of mutually orthogonal projections, we have $\varphi (p_0) \ldots \varphi (p_{n+1}) = 0$; this is equivalent to saying that $\bigcap_0^{n+1} \varphi (p_i)^{-1} ((0, \infty)) = \emptyset$.

\end{remarks}
\alten

\altbn{\label{cpr}}
We are now ready to define our notion of noncommutative covering dimension:

\begin{defn}
Let $A$ be a $C^*$-algebra. The completely positive rank of $A$ is less than or equal to $n$, $\cpr A \le n$, if the following holds: For $a_1 , \ldots , a_k \in A$ and $\varepsilon > 0$ there exists a c.p.\ approximation $(F,\psi , \varphi )$ for $\{a_1, \ldots , a_k \}$ within $\varepsilon$ such that the strict order of $\varphi$ does not exceed $n$.\\
We say $\cpr A = n$, if $n$ is the least integer such that $\cpr A \le n$.

\end{defn}
\alten

\altbn{\label{r-cpr}}
\begin{remark}
$\cpr A < \infty$ in particular implies the existence of c.p.\ approximations, thus $A$ must be nuclear.
\end{remark}
\alten

\altbn{\label{cpr-le-dim}}
The next result says that, for a commutative $C^*$-algebra, the completely positive rank is less than or equal to the covering dimension of the spectrum. We will see later that in fact they are equal.

\begin{prop}
Let $X$ be a second countable locally compact space. Then $\cpr (\Ch_0 (X)) \le \dim X$.
\end{prop}

\begin{proof}
 Let $n := \dim X$ and $a_1 , \ldots , a_k \in \Ch_0 (X)$ and $\varepsilon >0$ be given. Choose a compact subset $K \subset X$ such that $| a_i (x) | < \frac{\varepsilon}{3}$ for $x \in X \backslash K$ and $1 \le i \le k$. Then there is an open covering $U_0 , \ldots , U_r$ of $X$ with $U_0 = X \backslash K$ and such that $|a_i (x) - a_i (y) | < \frac{2}{3} \varepsilon \quad \mbox{for} \; x,y \in U_j \; , \; \forall  i,j$ (here we used that $K$ is compact).\\
Because $\dim X = n$, there is a refinement $(V_l)_{\{ 1 , \ldots , s \} }$ with the following properties:\\
(1) $|a_i (x) - a_i (y)| < \frac{2}{3} \varepsilon$ for $x,y \in V_l \; , \; \forall  i,l$,\\
(2) for every $\overline{l} \in \{ 1 , \ldots , s \}$ there exists $x_{\overline{l}} \in V_{\overline{l}}$ such that $x_{\overline{l}} \notin \bigcup_{l \neq \overline{l}} V_l$, \\
(3) $(V_l)_{ \{ 1 , \ldots , s \} }$ has strict order no greater that $n$.

Take a partition of unity $(h_l)_{ \{ 1 , \ldots , s \} }$ subordinate to $(V_l)_{\{ 1 , \ldots , s \} }$. Because $X$ is normal, there is a function $d : X \to [0,1]$ such that $d \in \Ch_0 (X)$ and $d (x) = 1$ for $x \in \{ x_1 , \ldots ,x_s \} \cup K$. \\
Set $h'_l := h_l \cdot d , \; l = 1 , \ldots , s$, and note that $h'_l \in \Ch_0 (X) , \, h'_l (x_l) = 1 \; \forall \, l$ and that $(h^{'-1}_l ((0,1]))_{ \{ 1 , \ldots , s \} }$ has strict order less than or equal to $n$. Also, $\sum^s_1 h'_l (x) = 1$ if $x \in K$. We are now prepared to define
\[
\begin{array}{l}
F := \C^s \\
\psi : \Ch_0 (X) \to F \, ,  \; \psi (a) := (a (x_1) , \ldots , a (x_s)) \; \mbox{and} \\
\varphi : F \to \Ch_0 (X) \, , \; \varphi (e_l) := h'_l \; .
\end{array}
\]
For $x \in X$ we obtain for every $i$:
\begin{eqnarray*}
  |\varphi \psi (a_i) (x) - a_i (x)| & \le & \Big| \sum^s_{l=1} a_i (x_l) h'_l (x) - \sum^s_{l=1} a_i (x) h'_l (x) \Big| \\
& & + \Big| \sum^s_{l=1} a_i (x) h'_l (x) - a_i (x) \Big| \\
& \le & \sum^s_{l=1} |a_i (x_l) - a_i (x)| h'_l (x) + |\sum h'_l (x) - 1| \, |a_i (x)| \\
& < & \frac{2}{3} \varepsilon \sum h'_l (x) + \frac{\varepsilon}{3} \\
& \le & \varepsilon \; ,
\end{eqnarray*}
because $|a_i (x_l) - a_i (x)| < \frac{2}{3} \varepsilon$, if $h'_l (x) > 0$ and $|a_i (x)| < \frac{\varepsilon}{3}$, if $|\sum h'_l (x)| < 1$. \\
Therefore, $\cpr (\Ch_0 (X)) \le n$ and the proof is complete.
\end{proof}
\alten

\altbn{\label{r-cpr-le-dim}}
\begin{remarks}
(i) It follows from the proof that, if $A$ is commutative, the approximating algebra $F$ may be chosen to be $\C^s$ for some $s \in \N$.\\
(ii) At the end of this section we will prove that in fact $\cpr (\Ch_0 (X)) = \dim X$ (Proposition \ref{dim=cpr}).
\end{remarks}
\alten 

\altbn{\label{matrices-le}} 
\begin{prop}
Let $X$ be a second countable locally compact space and $r \in \N$. Then $\cpr (\Ch_0 (X) \otimes M_r) \le \dim X$.
\end{prop}

\begin{proof}
  It obviously suffices to approximate elements of the form $a_j \otimes b_j$ for $a_j \in \Ch_0 (X)$, $b_j \in M_r$, $j=1, \ldots, k$. By Proposition \ref{cpr-le-dim} and Remark \ref{r-cpr-le-dim}(i) there is a c.p.\ approximation $(F,\psi, \varphi)$ within $\varepsilon$ for $\{a_1 , \ldots , a_k \} \subset \Ch_0 (X)$ with $n := \ord \varphi \le \dim X$ and $F \cong \C^s$ for some $s$. Then $(F \otimes M_r , \psi \otimes \id , \varphi \otimes \id)$ is a c.p.\ approximation for $\{a_1 \otimes b_1 , \ldots , a_k \otimes b_k\}$. We show $\ord (\varphi \otimes \id) \le \ord \varphi$:\\
We have $F \otimes M_r \cong M_r \oplus \ldots \oplus M_r$. Now let $\{ e_0 , \ldots , e_{n+1} \} \subset F$ be elementary, then each $e_i$ lives in a summand $M_r$, because it is minimal. This means $e_i = p_{e_i} \otimes q_i$, where $p_{e_i}$ is a generator of $\C^s$ and $q_i$ is a minimal projection in $M_r$.\\
If there are $i,j \in \{ 0 , \ldots , n+1 \}$ such that $p_{e_i} = p_{e_j}$, then $q_i \perp q_j$ because $e_i \perp e_j$, but then also $(\varphi \otimes \id) (e_i) \perp (\varphi \otimes \id) (e_j)$.\\
If the $p_{e_i}$ are pairwise orthogonal, then, because $\ord \varphi \le n$, there must be $i,j \in \{ 0 , \ldots , n+1 \}$ such that $\varphi (p_{e_i}) \perp \varphi (p_{e_j})$, thus $(\varphi \otimes \id) (e_i) \perp (\varphi \otimes \id) (e_j)$.\\
Therefore $\ord (\varphi \otimes \id) \le \ord \varphi$ and we are done.
\end{proof}
\alten

\altbn
\begin{nremark}
  Proposition \ref{matrices-le} will be generalized to continuous trace algebras in \cite{Wi}. There we will also show that in fact $\cpr (\Ch_0 (X) \otimes M_r) = \dim X$.
\end{nremark}
\alten

\altbn{\label{p-quotients}}
It follows immediately from the definition of covering dimension that, for a closed subset $K$ of some space $X$, $\dim K \le \dim X$. This carries over to the completely positive rank (note that quotients of unital $C^*$-algebras correspond to closed subsets of compact spaces): 

\begin{prop}
Let $A,B$ be nuclear $C^*$-algebras, $\pi : A \twoheadrightarrow B$ a surjective $*$-homomorphism. Then $\cpr B \le \cpr A$.
\end{prop}

\begin{proof}
  By the Choi--Effros lifting theorem \ref{lifting}, $\pi$ has a completely positive contractive lift $\sigma : B \to A$. \\
Let $b_1 , \ldots , b_k \in B$ and $\varepsilon > 0$ be given. Choose a c.p.\ approximation $(F , \psi , \varphi)$ for $\{ \sigma (b_1) , \ldots , \sigma (b_k)\}$ within $\varepsilon$ with $\ord \varphi \le \cpr A$. Then $(F , \psi \sigma , \pi \varphi)$ is a c.p.\ approximation of $\{ b_1 , \ldots , b_k \}$ within $\varepsilon$. $\ord (\pi \varphi) \le \ord \varphi$, because $\pi$ is a $*$-homomorphism, thus preserves orthogonality. 
\end{proof}
\alten

\altbn{\label{sums}}
\begin{prop}
If $A$ and $B$ are $C^*$-algebras, then
\[
\cpr (A \oplus B) = \max (\cpr A , \cpr B).
\]
\end{prop}
\begin{proof}
  Given $(a_1 , b_1) , \ldots , (a_k , b_k) \in A \oplus B$, $\varepsilon > 0$, choose c.p.\ approximations $(F_A , \psi_A , \varphi_A)$ and $(F_B , \psi_B , \varphi_B)$ within $\varepsilon$ for $\{ a_1 , \ldots , a_k \}$ and for $\{ b_1 , \ldots , b_k \}$, respectively, such that $\ord (\varphi_A) \le \cpr A$ and $\ord (\varphi_B) \le \cpr B$.\\
As a consequence, $(F_A \oplus F_B , \psi_A \oplus \psi_B , \varphi_A \oplus \varphi_B)$ is a c.p.\ approximation within $\varepsilon$ for $(a_1 ,b_1) , \ldots , (a_k ,  b_k)$. Since $\ord (\varphi_A \oplus \varphi_B) \le \max (\ord \varphi_A , \ord \varphi_B)$, we have $\cpr (A \oplus B) \le \max (\cpr A , \cpr B)$. Equality follows from Proposition \ref{p-quotients}.
\end{proof}
\alten

\altbn{\label{limits}}
The next result is an analogue of Proposition \ref{limits-dim}:

\begin{prop}
Let $(A_i , \varphi_{i , j})_{\N}$ be an inductive system of nuclear $C^*$-algebras and $A := \lim_{\to} A_i$.\\
Then $\cpr A \le \underline{\lim}\; \cpr A_i$.
\end{prop}

\begin{proof}
  Let $\varphi_i : A_i \to A$ be the induced maps, then $A := \overline{\bigcup_{\N} \varphi_i (A_i)}$. Given $a_1 , \ldots , a_k \in A$, $\varepsilon > 0$, choose $\bi \in \N$ such that there are $a'_1 , \ldots , a'_k \in \varphi_{\bi} (A_{\bi})$ with $\| a_j - a'_j \| < \frac{\varepsilon}{3}$, $j = 1 , \ldots , k$, and such that $\cpr A_{\bi} \le \underline{\lim} \; \cpr A_i$.\\
By Proposition \ref{p-quotients} there is a c.p.\ approximation $(F , \psi , \varphi)$ (of $\varphi_{\bi} (A_{\bi})$) within $\frac{\varepsilon}{3}$ for $\{ a'_1 , \ldots , a'_k \}$ with $\ord \varphi \le \cpr A_{\bi}$. After extending $\psi$ to all of $A$ by the Arveson extension theorem, $(F , \psi , \varphi)$ is a c.p.\ approximation within $\varepsilon$ for $\{ a_1 , \ldots , a_k \}$.
\end{proof}
\alten

\altbn{\label{subhomogeneous-limits}}
\begin{prop}
Let $A = \lim_{\longrightarrow} A_n$ be an inductive limit of $C^*$-algebras of the form 
\[
A_n = (M_{r (n,1)} \otimes \Ch_0 (X_{n,1})) \oplus \ldots \oplus (M_{r (n,k_n)} \otimes \Ch_0 (X_{n,k_n})) \, ,
\]
where the $X_{n,i}$ are second countable locally compact spaces.\\
Then $\cpr A \le \underline{\lim}_n \max_{i = 1 , \ldots , k_n} \dim (X_{n,i})$.
\end{prop}

\begin{proof}
  Let $\rho_n$ be the map from $A_n$ to $A$. It follows from Propositions \ref{p-quotients}, \ref{matrices-le} and \ref{sums} that $\cpr (\rho_n (A_n)) \le \max_{i = 1 , \ldots , k_n} (\dim (X_{n,i}))$. But now $A = \overline{\bigcup_{\N} \rho_n (A_n)}$, where the $\rho_n (A_n)$ are nested, so by Proposition \ref{limits} we have $\cpr A \le \underline{\lim} \; \cpr (\rho_n (A_n))$ and the assertion follows.
\end{proof}
\alten

\altbn{\label{examples}}
\begin{examples}
(i) Clearly, finite-dimensional $C^*$-algebras have completely positive rank zero. By Proposition \ref{limits} the same holds for $AF$ algebras. We will see that in fact $\cpr A =0$ iff $A$ is $AF$ (Theorem \ref{AF}).\\
(ii) From Proposition \ref{subhomogeneous-limits} it follows that $AT$ algebras (i.e.\ inductive limits of direct sums of matrix algebras over $S^1$) have completely positive rank less than or equal to one. In particular, if $A$ is a Bunce--Deddens algebra (cf.\ \cite{BD}) or an irrational rotation algebra (cf.\ \cite{Ri2} and \cite{EE}), $\cpr A =1$, since these are $AT$ but not $AF$.\\
(iii) In \cite{Bl2} Blackadar constructed a simple, unital, projectionless $C^*$-algebra, which is easily seen to have completely positive rank one, although it is not $AT$.
\end{examples}
\alten

The rest of this section is devoted to the proof of the fact that in Proposition \ref{cpr-le-dim} we even have equality (Proposition \ref{dim=cpr}). However, this turns out to be not so easy, because we have to construct refinements of open coverings from c.p.\ approximations $(F , \psi , \varphi)$ with not necessarily commutative algebras $F$ (if $F \cong \C^k$, this yields an open covering of $X$ by associating to each generator of $\C^k$ an open subset of $X$). If $F$ is the direct sum of matrix algebras, it is not clear how to obtain the right open covering of $X$ from $(F , \psi , \varphi)$ (an arbitrary copy of $\C^k$ in $F$ will not be good enough, neither will be the copy generated by the central projections of $F$). The key observation is the next lemma. It implies that the size of the matrix algebras in $F$ cannot be too large.

\altbn{\label{dichotomy}}
\begin{lemma}
Let $A$ be a $C^*$-algebra and $\varphi : M_r \to A$ be a c.p.\ contraction. Then either $\ord \varphi =0$ or $\ord \varphi = r - 1$.
\end{lemma}

\begin{proof}
If $r=1$, then $\ord \varphi =0$ and there is nothing to show, so we assume $r>1$. Let $n:= \ord \varphi$ and suppose $n \neq 0$. We (inductively) construct an elementary set $E_r \subset M_r, \; |E_r| =r$, such that $\varphi (e) \, \varphi (\overline{e}) \neq 0$ for $e, \, \overline{e} \in E_r$. This will imply $n> r-2$; we obviously have $n< r-1$, so $n = r-1$ if $n \neq 0$.\\
There are orthogonal minimal projections $e_1 , \, e_2 \in M_r$ such that $\varphi(e_1) \, \varphi(e_2) \neq 0$; set $E_2 := \{ e_1, \, e_2\}$. Next suppose that for some $k \in \{ 2, \ldots , r-1 \}$ we have constructed an elementary set $E_k = \{ e'_1, \ldots e'_k \} \subset M_r$ such that
\[
\varphi (e'_i) \, \varphi (e'_j) \neq 0 \; \forall i, \, j \in \{1, \ldots , k\}.
\]
Choose a minimal projection $e'_{k+1} \in M_r$ orthogonal to $E_k$; this is possible since $k<r$. Take $e'_k \in E_k$ and identify $M_2$ with $(e'_k + e'_{k+1}) M_r (e'_k + e'_{k+1})$. There is an open neighborhood $\Uh$ of $\be_{M_2}$ in $U (M_2)$ such that
\[
\varphi (u^* e'_k u) \, \varphi (e'_i) \neq 0 \; \forall i \in \{1, \ldots , k-1 \}, \; u \in \Uh \,.
\]
Note that $\{ e'_1 , \ldots e'_{k-1}, \, u^* e'_k u, \, u^* e'_{k+1} u \} \subset M_r$ is elementary for all $u \in \Uh$. Set $\Fh_i := \{ u\in \Uh \, | \varphi (u^* e'_{k+1} u) \perp \varphi (e'_i) \}$ for $i = 1, \ldots, k-1$ and $\Fh := \bigcup^{k-1}_{i=1} \Fh_i$, then the $\Fh_i$ and $\Fh$ are closed in $\Uh$. If $\Uh \backslash \Fh$ was empty, this would mean
\[
\varphi (u^* e'_{k+1} u) \perp \varphi(e'_1) \; \forall u \in \Uh \, ;
\]
but then by Lemma \ref{unitary-neighborhoods} and Remark \ref{dominated-orthogonality} we would have $\varphi (e'_k + e'_{k+1}) \perp \varphi (e'_1)$, a contradiction to $ \varphi (e'_k) \, \varphi (e'_1) \neq 0$.\\ 
So let $\Vh_1 := \Uh \backslash \Fh_1$; just as above we obtain that $\Vh_2 := \Vh_1 \backslash \Fh_2$ must be nonempty. Inductively we get that $\Vh := \Uh \backslash \Fh \neq \emptyset$. $\Vh$ is open in $U(M_2)$ and we have
\[
\varphi (u^* e'_{k+1} u) \, \varphi (e'_i) \neq 0 \; \forall u \in \Vh, \, i \in \{1, \ldots k-1\} \,.
\]
Choose some $\bar{u} \in \Vh$ and set $e''_j := \bar{u}^* e'_j \bar{u}, \;j = k, \, k+1$. It might still happen that $\varphi (e''_k) \perp \varphi (e''_{k+1})$. To fix this, identify $M_2$ with $(e'_{k-1} + e''_k) M_r (e'_{k-1} + e''_k)$.\\ 
There is an open neighborhood $\Uh' \subset U(M_2)$ of $\be_{M_2}$ such that $\varphi (e'_i) \, \varphi (u^* e'_{k-1} u)$, $ \varphi (e''_{k+1}) \, \varphi (u^* e'_{k-1} u)$, $ \varphi (u^* e''_k u) \, \varphi(u^* e''_{k-1} u)$ and $\varphi (e'_i) \, \varphi(u^* e''_k u)$ are all nonzero for $i = 1, \ldots , k-2$ and $u \in \Uh'$.\\
Now there must be a nonempty open subset $\Vh' \subset \Uh'$ such that
\[
\varphi (u^* e''_k u) \, \varphi (e''_{k+1}) \neq 0 \; \forall u \in \Vh' \, ,
\]
for otherwise again by Lemma \ref{unitary-neighborhoods} und Remark \ref{dominated-orthogonality} we would have $\varphi ( e'_{k-1}) \perp \varphi (e''_{k+1})$, contradicting our construction. Choose some $\bar{v} \in \Vh'$ and define
\[
E_{k+1} := \{ e_1' , \ldots , e_{k-2}', \bar{v}^* e_{k-1}' \bar{v} , \bar{v}^* e_k'' \bar{v} , e_{k+1}'' \} \, ,
\]
then $E_{k+1}$ is elementary by construction, $| E_{k+1} | = k+1$ and we have
\[
\varphi (e) \, \varphi (\bar{e}) \neq 0 \; \forall e, \bar{e} \in E_{k+1} \, .
\]
Induction then yields an elementary set $E_r$ with the desired properties.
\end{proof}
\alten

\altbn{\label{r-dichotomy}}
\begin{remarks}
(i) The lemma might also be read as follows: Let $\varphi : M_r \to A$ be c.p.c.\ with $\ord \varphi \le n$. Then $\ord \varphi = 0$ or $r \le n+1$.\\
(ii) If $A$ is commutative, then we always have $r \le n+1$, since $\ord \varphi =0$ if and only if $r = 1$.
\end{remarks} 
\alten 

\altbn{\label{matrixrank}}
\begin{lemma}
Let $a_1 , \ldots , a_k \in M_n$ be positive elements with $\sum a_i \le \be$ and $\| a_i \| > \frac{n}{n+1}$. Then $k \le n$.
\end{lemma}

\begin{proof}
  $1 = \tr (\be) \ge \tr \left( \sum_i a_i \right) = \sum_i \tr (a_i) \ge \frac{1}{n} \sum_i \| a_i \| > \frac{1}{n} \cdot k \cdot \frac{n}{n+1} = \frac{k}{n+1}$, so $k < n+1$.
\end{proof}
\alten

\altbn{\label{orthogonal-projections2}}
\begin{lemma}
Let $K$ be a natural number and $\beta > 0$. Then there is $\alpha > 1$, such that the following holds: \\
Given projections $q_1 , \ldots , q_k , \, k \le K$, in a $C^*$-algebra $A$ with $\| \sum^k_1 q_i\| \le \alpha$, there exist pairwise orthogonal projections $p_1 , \ldots , p_k \in A$ satisfying 
\[
\|p_i - q_i\| \le \beta , \; i = 1 , \ldots , k.
\]
 If $\| \sum^k_1 q_i\| \le 1$ , then the $q_i$ already are pairwise orthogonal themselves.
\end{lemma}

\begin{proof}
  For any $\alpha \ge 1$ and for $i \in \{ 1 , \ldots , k \}$ we have
\[
q_i \left( \sum^{i-1}_1 q_l \right) q_i \le q_i \left( \sum^k_1 q_l \right) q_i - q_i \le (\alpha -1) q_i \; ,
\]
therefore
\[
\Big\| q_i \sum^{i-1}_1 q_l \Big\| \le \Big\| q_i \left( \sum^{i-1}_1 q_l \right)^{\halb} \Big\| \; \Big\| \left( \sum^{i-1}_1 q_l \right)^{\halb} \Big\| \le (\alpha -1)^{\halb} \alpha^{\halb} =: \delta \; .
\]
Now inductively construct pairwise orthogonal projections $p_i$ and $\delta_i \ge 0$, such that
\[
\| p_i - q_i \| \le \delta_i \; , \; i = 1 , \ldots , k \; .
\]
First define $p_1 := q_1 \; , \; \delta_1 := 0$. Next suppose $p_l$, $\delta_l , \, l = 1, \ldots , i-1$ already are constructed. We then have
\begin{eqnarray*}
  \Big\| q_i \left( \sum^{i-1}_1 p_l \right) \Big\| & \le & \Big\| q_i \sum^{i-1}_1 (p_l - q_l) \Big\| + \Big\| q_i \sum^{i-1}_1 q_l \Big\| \\
& \le & \sum^{i-1}_1 \delta_l + \delta \; .
\end{eqnarray*}
Corollary \ref{orthogonal-projections} now yields a projection $p_i \perp \sum^{i-1}_1 p_l$ with 
\[
\| p_i - q_i \| \le 14 \left( \sum^{i-1}_1 \delta_l + \delta \right) := \delta_i.
\]
Obviously $\delta_i \le \beta $ for $ i = 1 , \ldots , k$, if $\alpha$ has been chosen small enough.
\end{proof}
\alten

\altbn{\label{dim=cpr}}
\begin{prop}
Let $X$ be locally compact and second countable. Then $\dim X = \cpr (\Ch_0 (X))$.
\end{prop}

\begin{proof}
  $\cpr (\Ch_0 (X)) \le \dim X$ by Proposition \ref{cpr-le-dim}, so we have to show $\dim X \le \cpr (\Ch_0 (X))$. Unfortunately, the proof is a little technical, so we briefly sketch the strategy:\\
Given an open covering $(U_{\gamma})_{\Gamma}$ of $X$, choose a refinement $(V_{\lambda})_{\Lambda}$, where the $V_{\lambda}$ are suffiently small. Take a partition of unity $(h_{\lambda})_{\Lambda}$ subordinate to $(V_{\lambda})_{\Lambda}$ and a c.p.\ approximation $(F , \psi , \varphi)$ with $\ord \varphi \le \cpr (\Ch (X))$ for the $h_{\lambda}$. By Lemma \ref{dichotomy}, $F$ must be a direct sum of matrices of rank no larger than $\cpr (\Ch (X)) + 1$. On suitable subsets of $\Lambda$ one now can introduce equivalence relations to combine the $V_{\lambda}$ to a new open covering $(\tilde{V}^{(i)}_j)$, which still refines $(U_{\gamma})_{\Gamma}$. One can use Lemma \ref{orthogonal-projections2} and the special properties of $(F, \psi , \varphi)$ to produce a refinement $(W^{(i)}_j)$ of $(\tilde{V}^{(i)}_j)$, the order of which is less than or equal to $\ord \varphi$. Then $(W^{(i)}_j)$ will be the desired refinement of $(U_{\gamma})_{\Gamma}$. 

$X$ is second countable, thus $X = \bigcup_{\N} A_m$ for compact subsets $A_m \subset X$. By Proposition \ref{p-quotients} we have $\cpr (\Ch (A_m)) \le \cpr (\Ch_0 (X)) \; \forall m$. Now if we can show that $\dim A_m \le \cpr (\Ch (A_m))$ we are done, for the countable sum theorem \ref{countable-sums} then yields $\dim X \le \overline{\lim} \dim (A_m) \le \cpr (\Ch_0 (X))$.\\ 
Therefore we may assume $X$ to be compact (and second countable); in particular we may choose a metric $\mathrm{d}$ on $X$. Set $n := \cpr (\Ch (X))$. 

Choose a partition of unity $(f_{\gamma})_{\Gamma}$ subordinate to the given (finite) open covering $(U_{\gamma})_{\Gamma}$ of $X$. Set $\varepsilon := \frac{1}{|\Gamma|}$, then there exists a $\delta > 0$, such that $| f_{\gamma} (x) - f_{\gamma} (y) | < \varepsilon \; \forall \gamma \in \Gamma$ and $ x,y \in X$ with $\mathrm{d} (x,y) < \delta$.\\
Again because $X$ is compact, there exists a finite open covering $(V_{\lambda})_{\Lambda}$ with 
\[
\diam (V_{\lambda}) < \frac{\delta}{3(n+1)} \; \forall \lambda \in \Lambda;
\]
 let $(h_{\lambda})_{\Lambda}$ be a partition of unity subordinate to $(V_{\lambda})_{\Lambda}$.

We have to define some constants:\\
First set $C := \frac{1}{2 (n+1)} , \; \beta := \frac{C}{2} = \frac{1}{4(n+1)}$. Then choose $1 < \alpha \le \frac{n+2}{n}$ as in Lemma \ref{orthogonal-projections2} (with $K = n+1$) and define $\vartheta := \frac{1}{\alpha} , \; \eta := (1 - \vartheta) \frac{C}{2}$. Note that $\frac{\eta}{C} = \frac{1-\vartheta}{2} \le \halb \left( 1 - \frac{n}{n+2} \right) = \frac{1}{n+2}$. Take a c.p.\ approximation $(F , \psi, \varphi) , \; \ord \varphi \le n$, within $\eta$ for $\{ \sum_{\lambda \in \Lambda'} h_{\lambda} \tei \Lambda' \subset \Lambda \}$. Then $F \cong F_1 \oplus \ldots \oplus F_m , \; F_j = M_{r_j}$, and by Lemma \ref{dichotomy} we have $r_j \le n+1 , \; j = 1 , \ldots , m$. Let $\be_j$ be the unit of $F_j$ and $\psi_j$ be the $j$'th component of $\psi$. Now define
\[
A_j := \{ x \in X \tei \varphi (\be_j) (x) > C \}
\]
and
\[
\Lambda_j := \{ \lambda \in \Lambda \tei V_{\lambda} \cap A_j \neq \emptyset \} \; .
\]
Denote by $\sim_j$ the equivalence relation on $\Lambda_j$ generated by the relation 
\[
\lambda_1 \sim_j \lambda_2 \iff V_{\lambda_1} \cap V_{\lambda_2} \cap A_j \neq \emptyset.
\] 

Let $\Lambda^{(1)}_j , \ldots , \Lambda^{(s_j)}_j \subset \Lambda_j$ be the equivalence classes with respect to $\sim_j$, then $\Lambda_j = \dot{\bigcup} \Lambda^{(i)}_j$.\\
Set $\tilde{V}^{(i)}_j := \bigcup_{\lambda \in \Lambda^{(i)}_j} V_{\lambda} \cap A_j , \; i = 1 , \ldots , s_j$; then $A_j \subset \bigcup_i \tilde{V}^{(i)}_j$ and $\tilde{V}^{(i_1)}_j \cap \tilde{V}^{(i_2)}_j = \emptyset$, if $i_1 \neq i_2$.

Setting $h^{(i)}_j := \sum_{\lambda \in \Lambda^{(i)}_j} h_{\lambda}$, we have $\| \varphi \psi (h^{(i)}_j) - h^{(i)}_j \| < \eta$ and $h^{(i)}_j (x) = 1$ for $x \in \tilde{V}^{(i)}_j$. It follows that
  \begin{eqnarray*}
    \varphi (\be_j - \psi_j (h^{(i)}_j)) (x) & \le & \varphi (\be_F - \psi (h^{(i)}_j)) (x) \\
& \le & (\be_{\Ch_0 (X)} - \varphi \psi (h^{(i)}_j)) (x) \\
& < & 1 - h^{(i)}_j (x) + \eta = \eta \; .
  \end{eqnarray*}
Define $q^{(i)}_j := g_{\vartheta} (\psi_j (h^{(i)}_j)) \in F_j$ (cf.\ \ref{functions}); $q^{(i)}_j$ is a projection and for $x \in \tilde{V}^{(i)}_j$ we obtain:
\begin{eqnarray}
  \vartheta \cdot \varphi (\be_j - q^{(i)}_j) (x) + \varphi (q^{(i)}_j) (x) & \ge & \varphi (\psi_j (h^{(i)}_j)) (x) > \varphi (\be_j) (x) - \eta \nonumber \\
\Rightarrow (1 - \vartheta) \cdot \varphi (\be_j - q^{(i)}_j) (x) & < & \eta \nonumber \\
\Rightarrow \hspace*{1.3cm} \varphi (\be_j - q^{(i)}_j) (x) & < & \frac{\eta}{1-\vartheta} = \frac{C}{2} \; \; \forall\, x \in \tilde{V}^{(i)}_j \; . \label{1}
\end{eqnarray}
We have
\[
\sum^{s_j}_{i=1} q^{(i)}_j \le \frac{1}{\vartheta} \cdot \sum^{s_j}_{i=1} \psi_j (h^{(i)}_j) \le \frac{1}{\vartheta} \cdot \be_j = \alpha \cdot \be_j \; ,
\]
so by Lemma \ref{orthogonal-projections2} (and by the choice of $\alpha$) there are pairwise orthogonal projections $p^{(1)}_j , \ldots , p^{(s_j)}_j \in F_j$ with $\| p^{(i)}_j - q^{(i)}_j \| < \beta , \; i = 1 , \ldots , s_j$. 

We are now prepared to define
\[
W^{(i)}_j := \{ x \in X \tei \varphi (p^{(i)}_j) (x) > C \} \subset A_j \; ;
\]
we will show that $(W^{(i)}_j)_{j = 1 , \ldots , m \atop i = 1 , \ldots , s_j}$ is an open covering that refines $(U_{\gamma})_{\Gamma}$. Clearly, the $W^{(i)}_j$ are open. By Remark \ref{strict-order} (ii), the order of $(W^{(i)}_j)$ is less than or equal to $n$, because $\ord \varphi \le n$.

$(W^{(i)}_j)$ covers $X$:\\
Let $x$ be in $X$. $\be_F$ can be decomposed into sums of orthogonal projections:
\[
\be_F = \sum_{j = 1 , \ldots , m \atop i = 1 , \ldots , s_j} p^{(i)}_j + \sum^m_{j=1} \left( \be_j - \sum^{s_j}_{i=1} p^{(i)}_j \right) \; ;
\]
since $\ord \varphi \le n$, for at most $n+1$ projections of the form $p^{(i)}_j$ or $(\be_j - \sum_i p^{(i)}_j)$, $\varphi (p^{(i)}_j) (x)$ or $\varphi \left( \be_j - \sum_i p^{(i)}_j \right) (x)$ can be nonzero (again we used Remark \ref{strict-order} (ii)).

Now for all $i ,\, j$ we have $\varphi (\be_j - \sum_i p^{(i)}_j) (x) \le C$: this is clear if $x \notin A_j$; if $x \in A_j$, then $x \in \tilde{V}^{(i')}_j$ for some $i'$, but then (\ref{1}) yields
\[
\varphi \left( \be_j - \sum_i p^{(i)}_j \right) (x) \le \varphi (\be_j - p^{(i')}_j ) (x) \stackrel{(*)}{<} \frac{\eta}{1-\vartheta} + \beta = C < \frac{1-\eta}{n+1} \; .
\]
We also have $| \varphi (\be_F) (x) - 1 | < \eta$, thus for at least one $p^{(i)}_j$ we obtain 
\[
\varphi (p^{(i)}_j) (x) \ge \frac{1-\eta}{n+1} > C \, ,
\]
i.e.\ $x \in W^{(i)}_j$ and $(W^{(i)}_j)$ is a covering.

It remains to check that $(W^{(i)}_j)$ is a refinement of $(U_{\gamma})_{\Gamma}$. We first show $W^{(i)}_j \subset \tilde{V}^{(i)}_j \; \forall i,j$:\\
Be $x \in W^{(i)}_j$. Then $x \in A_j$, thus $x \in \tilde{V}^{(i')}_j$ for some $i'$. If $i \neq i'$, then again (\ref{1}) gives
\[
C < \varphi (p^{(i)}_j) (x) \le \varphi (\be_j - p^{(i')}_j) (x) \stackrel{(*)}{<} \frac{C}{2} + \frac{C}{2} = C\, ,
\]
a contradiction, so $i = i'$ and $x \in \tilde{V}^{(i)}_j$. 

Next we show $\mathrm{d} (x,y) < \delta$ for $x,y \in \tilde{V}^{(i)}_j$:\\
Suppose there are $x,y \in \tilde{V}^{(i)}_j$ with $\mathrm{d} (x,y) \ge \delta$. Recall that 
\[
\tilde{V}^{(i)}_j = \bigcup_{\lambda \in \Lambda^{(i)}_j} V_{\lambda} \cap A_j
\]
and that $\Lambda^{(i)}_i \subset \Lambda$ is an equivalence class with respect to $\sim_j$. This means there are $\lambda_1 , \ldots , \lambda_t \in \Lambda^{(i)}_j$ such that $x \in V_{\lambda_1} , \; y \in V_{\lambda_t}$ and 
\[
V_{\lambda_l} \cap V_{\lambda_{l+1} } \cap A_j \neq \emptyset , \; l = 1 , \ldots , t-1.
\]
But $\diam (V_{\lambda}) < \frac{\delta}{3(n+1)} \; \forall \lambda$, thus for $k = 1 , \ldots , n$ one can choose $x_k \in \tilde{V}^{(i)}_j$, such that 
\[
\frac{k\delta}{n+1} - \frac{\delta}{3(n+1)} \le \mathrm{d} (x_1 , x_k) \le \frac{k\delta}{n+1}.
\]
Setting $x_0 := x , \; x_{n+1} := y$, one obtains $x_0 , \ldots , x_{n+1} \in \tilde{V}^{(i)}_j$ with $\mathrm{d} (x_{k_1} , x_{k_2}) \ge \frac{2\delta}{3 (n+1)}$ if $k_1 \neq k_2$. Set
\[
\Lambda_k := \{ \lambda \in\Lambda \tei x_k \in V_{\lambda} \} , \; \tilde{h}_k := \sum_{\lambda \in \Lambda_k} h_{\lambda} , \; k = 0 , \ldots , n+1 \; ,
\]
then $\Lambda_{k_1} \cap \Lambda_{k_2} = \emptyset$ if $k_1 \neq k_2$, because $\mathrm{d} (x_{k_1} , x_{k_2}) \ge \frac{2\delta}{3(n+1)}$ and $\diam V_{\lambda} < \frac{\delta}{3(n+1)}$, so $\sum_k \tilde{h}_k \le \sum_{\Lambda} h_{\lambda}$. Also, $\tilde{h}_k (x_k) = 1$ and $\| \varphi \psi (\tilde{h}_k) - \tilde{h}_k \| < \eta$, thus 
\[
0 \le \varphi (\be_j) (x_k) - \varphi (\psi_j (\tilde{h}_k)) (x_k) < \eta.
\]
But $x_k \in A_j$, therefore $\varphi (\be_j) (x_k) > C$; altogether we obtain:
\[
\begin{array}{ll}
 & \varphi (\be_j) (x_k) - \eta < \varphi (\psi_j (\tilde{h}_k)) (x_k) \le \| \psi_j (\tilde{h}_k)\| \cdot \varphi (\be_j) (x_k) \\
\\
\Rightarrow & 1 - \| \psi_j (\tilde{h}_k) \| < \frac{\eta}{\varphi (\be_j) (x)} < \frac{\eta}{C} \le \frac{1}{n+2} \, , \; k = 0 , \ldots , n+1 \; .
\end{array}
\]
We thus have $\| \psi_j (\tilde{h}_k)\| > \frac{n+1}{n+2}$, furthermore 
\[
\sum^{n+1}_{k=0} \psi_j (\tilde{h}_k) \le \sum_{\lambda\in \Lambda} \psi_j (h_{\lambda}) \le \be_j.
\]
But $\psi_j (\tilde{h}_k) \in F_j \cong M_{r_j} , \; k = 0 , \ldots , n+1$, and $r_j \le n+1$, so $\frac{n+1}{n+2} \ge \frac{r_j}{r_j+1}$. Now we can apply Lemma \ref{matrixrank} to obtain $n+2 \le r_j$. On the other hand, $\ord \varphi \le n$, so in particular $\ord (\varphi |_{M_{r_j}}) \le n$. Then by Remark \ref{r-dichotomy} we have $r_j \le n+1$, a contradiction. Therefore $\mathrm{d} (x,y) < \delta$. \\
Finally it turns out that $\tilde{V}^{(i)}_j \subset U_{\gamma}$ for some $\gamma \in \Gamma$: Take $x \in \tilde{V}^{(i)}_j$, then there is $\gamma \in \Gamma$ such that $f_{\gamma} (x) \ge \frac{1}{|\Gamma|}$. Then for any $y \in \tilde{V}^{(i)}_j$ we have $\mathrm{d} (x,y) < \delta$, therefore $|f_{\gamma} (y) - f_{\gamma} (x)| < \varepsilon$, and $f_{\gamma} (y) > \frac{1}{|\Gamma|} - \varepsilon = 0$, so $y \in U_{\gamma}$. But this means $\tilde{V}^{(i)}_j \subset U_{\gamma}$.

Thus we have shown that $(W^{(i)}_j)_{j = 1 , \ldots , m \atop i =1 , \ldots , s_j}$ covers $X$, is of order less than or equal to $n$ and refines $(U_{\gamma})_{\Gamma}$. Therefore, $\dim X \le n$. 
\end{proof}
\alten
\section[{\sc The zero-dimensional case}]{\large \sc The zero-dimensional case}

\subsection[{\rm Maps of strict order zero}]{\sc Maps of strict order zero}

Of course it is an important question what it means for a map to be of strict order $n$. Below we explicitly describe maps of strict order zero.

\bn{\label{ord=0}}
\begin{prop}
Let $A,F$ be $C^*$-algebras, $F = M_{r_1} \oplus \ldots \oplus M_{r_k}$ finite-dimensional and $\varphi : F \to A$ a completely positive contraction of strict order zero.\\
a) There are closed subspaces $X_1 , \ldots , X_k \subset (0,1]$ and an isomorphism $\alpha : C^* (\varphi (F)) \cong \bigoplus^k_1 \Ch_0 (X_i, M_{r_i})$ such that 
\[
\alpha \verk \varphi (y_i) (t_i) = t_i \cdot y_i 
\]
for $t_i \in X_i$ and $y_i \in M_{r_i} \subset F$.\\
b) If $\varphi (\be_F)$ is a projection, then $\varphi$ is a $*$-homomorphism.\\
c) For $\delta > 0$ there is $\gamma > 0$ such that the following holds:\\
If $\| \varphi (\be_F) - \varphi(\be_F)^2 \| < \gamma$, then there is a $*$-homomorphism $\varphi' : \, F \to A$ such that $\| \varphi' - \varphi \| < \delta$.
\end{prop}
\en

\bn{\label{r-ord=0}}
\begin{remark}
With the identification $\Ch_0(X_i, M_{r_i}) \cong \Ch_0(X_i) \otimes M_{r_i}$, in a) we could also write 
\[
\alpha \verk \varphi = {\te \bigoplus\limits^k_1} \, h_i \otimes \id_{M_{r_i}} \, ,
\]
where $h_i$ denotes the identity map on $X_i$. 
\end{remark}
\en

\begin{proof} (of \ref{ord=0})\\
a) (i)  Since $\ord \varphi = 0$, we have $\varphi (M_{r_i}) \perp \varphi (M_{r_j})$ for $i \neq j$, thus $C^* (\varphi (F)) = \bigoplus^k_1 C^* (\varphi (M_{r_i}))$ and we may assume $F = M_r$ for some $r \in \N$. Further we may assume that $A = C^* (\varphi (M_r))$ and that $A \subset \Bh (H)$ acts nondegenerately on some Hilbert space $H$. Then in particular the support projection of $h := \varphi (\be_{M_r})$ is $\be_H$. (If $\varphi$ is the zero map, then the closed subspace $X$ of $(0,1]$ will just be the empty set and there is nothing to show.)  

(ii) $[h , \varphi (x)] = 0 \; \forall x \in M_r$: It suffices to consider $x \ge 0$. But then $x = \lambda_1 \cdot e_1 + \ldots + \lambda_r \cdot e_r$ for some elementary set $\{ e_1 , \ldots , e_r \} \subset M_r$ and $0 \le \lambda_1 , \ldots , \lambda_r \le \| x \|$.\\
Now $h = \varphi (\be_{M_r}) = \sum^r_1 \varphi (e_i)$, and since $\varphi (e_i) \varphi (e_j) = 0$ for $i \neq j$, we see that $h \varphi (x) = \varphi (x) h = \sum \lambda_i \varphi (e_i)^2$.

(iii) Lemma \ref{unital-compression} says that $\varphi (\, .\,) = h^{\halb} \sigma (\, .\,) h^{\halb}$ for some u.c.p.\ map $\sigma : M_r \to \Bh (H)$, but a glance at the proof tells us that in fact
\[
\left( h + \frac{1}{n} \cdot \be_H \right)^{-\halb} \varphi (x) \left( h + \frac{1}{n} \cdot \be_H \right)^{-\halb} \stackrel{\rm s.o.}{\longrightarrow} \sigma (x) \; \forall x \in M_r \; .
\]
We obviously have $[h , \sigma (x)] = 0 \; \forall x \in M_r$.\\
Note that $\ord \sigma = 0$: Take rank-one projections $e \perp f \in M_r$. Since $\supp h = \be_H$, we see that 
\[
\sigma (e) \sigma (f) = 0 \iff h \sigma (e) \sigma (f) h = 0,
\]
but 
\[ 
h \sigma (e) \sigma (f) h = h^{\halb} \sigma (e) h^{\halb} h^{\halb} \sigma (f) h^{\halb} = \varphi (e) \varphi (f) = 0,
\]
because $\ord \varphi = 0$.

Using the universal property of the maximal tensor product and the fact that $C^* (h)$ and $C^* (\sigma (M_r))$ commute, it is straightforward to check that $C^* (\varphi (M_r)) \cong C^* (h) \otimes C^* (\sigma (M_r))$ and that, under this identification, $\varphi (x) = h \otimes \sigma (x)$.\\
Note that $C^* (h)$ is a quotient of $\Ch_0 ((0,1])$, the universal $C^*$-algebra generated by a positive element of norm $\le 1$; so $C^* (h) \cong \Ch_0 (X)$ for some closed $X \subset (0,1]$, and the isomorphism is given by sending $h$ to $\id_X$. 

(iv) The only thing left to show is that $\sigma$ is a $*$-homomorphism, for then $C^* (\sigma (M_r)) = \sigma (M_r) \cong M_r$ and with all these identifications $\varphi$ is given by $x \mapsto \id_X \otimes x$.

For any rank-one projection $e \in M_r$ we have $\sigma (e) \perp \sigma (\be_{M_r} - e)$ because $\ord \sigma = 0$; now since $\sigma (\be_{M_r}) = \be_H$ is a projection, so is $\sigma (e)$. But then by Proposition \ref{multiplicativity}, $\sigma (ex) = \sigma (e) \sigma (x) \; \forall x \in M_r$; since $e$ was arbitrary and $M_r$ is spanned by its minimal projections, $\sigma$ is multiplicative on all of $M_r$.

b) If $\varphi(\be_F)$ is a projection, then so is $\varphi (\be_{M_{r_i}})$ for each $i$ (the $\varphi (\be{M_{r_i}})$ are mutually orthogonal and add up to $\varphi (\be_F)$). But then we have for each $t_i \in X_i$
\[
t_i \cdot \be_{M_{r_i}} \stackrel{\rm a)}{=} \alpha \verk \varphi (\be_{M_{r_i}}) (t_i) = \be_{M_{r_i}} \, ,
\]
which implies that $t_i = 1$ and $X_i = \{ 1 \}$ for all $i$ for which $X_i \neq \emptyset$, and this in turn means that $\varphi$ is a $*$-homomorphism.

c) We will obtain $\varphi'$ from $\varphi$ by setting $\varphi' ( \, . \,) := c \varphi ( \, . \,) c$ for a suitable element $c \in C^* (\varphi (F)) = \bigoplus_1^k C^* (\varphi (M_{r_i}))$. Then $\| \varphi' - \varphi \| = \max_{i=1, \ldots, k} \| \varphi'_i - \varphi_i \|$, where $\varphi_i = \varphi |_{M_{r_i}}$ and $\varphi'_i = \varphi'|_{M_{r_i}}$. Note also that $\varphi'_i (\, . \,) = c_i \varphi_i (\, . \,) c_i$, where $c_i \in C^* (\varphi (M_{r_i}))$ is the i-th component of $c$. We may therefore again assume $F = M_r$ for some $r \in \N$.

Now suppose $\| \varphi ( \be_{M_r}) - \varphi ( \be_{M_r})^2 \| < \gamma$ (it will soon become clear how small $\gamma$ has to be). By Proposition \ref{almost-projections} there exists a projection $p \in C^* (\varphi ( \be_{M_r})) \subset A$ satisfying $\| p - \varphi ( \be_{M_r}) \| < 2 \gamma$. For $c := (p \varphi ( \be_{M_r}) p)^{- 1/2} \in C^*( \varphi( \be_{M_r}))_+ $ we have $p = c \varphi (\be_{M_r}) c$ and $\| p - c \| < 4 \gamma$.

Let $d \in (p A p)_+$ be another element with $d^2 = (p \varphi (\be_{M_r}) p)^{-1}$, then by uniqueness of the inverse $d^2 = c^2$ and $d = (d^2)^{1/2} = (c^2)^{1/2} =c$.

Now define $\varphi' : M_r \to pAp$ as above by $\varphi' (x) := c \varphi (x) c$. Then $\varphi' (\be_{M_r}) = p$ and $\| \varphi' - \varphi \|$ is small, because $\| c - \varphi (\be_{M_r}) \| < 6 \gamma$ and because for $x \in (M_r)_+ , \; \| x \| \le 1$, we have
\begin{eqnarray}
\label{2}
  \| \varphi (x) - \varphi (\be_{M_r}) \varphi (x) \|^2 & = & \| (\be - \varphi (\be_{M_r})) \varphi (x)^2 (\be - \varphi (\be_{M_r})) \| \nonumber \\ 
 & \le & \| (\be - \varphi (\be_{M_r})) \varphi (\be_{M_r}) (\be - \varphi (\be_{M_r})) \| \nonumber \\
 & < & \gamma \; .
\end{eqnarray}
Therefore
\begin{eqnarray*}
\| \varphi' (x) - \varphi (x)\| & \le & \| c \varphi (x) c - \varphi (\be_{M_r}) \varphi (x) \varphi (\be_{M_r}) \| \\
& & + \| \varphi (\be_{M_r}) \varphi (x) \varphi (\be_{M_r}) - \varphi (x) \| \\
& < & 12 \gamma + 2 \gamma^{\halb} \; .
\end{eqnarray*}
Thus, if only $\gamma$ is small enough, we have $\| \varphi' - \varphi \| < \delta$. 

Let $\{ e_1 , \ldots , e_{r} \} \subset M_r$ be some elementary set, then $\sum e_i = \be_{M_r}$ and $\varphi (e_i) \perp \varphi (e_j)$, $i \neq j$, and for every positive continuous function $f$ one has $f (\varphi (e_i)) \perp f (\varphi (e_j))$ , $i \neq j$ and $f (\varphi (\be_{M_r})) = \sum f (\varphi (e_i))$.

Using (\ref{2}) one obtains
\begin{eqnarray*}
  \| \varphi (e_j) - \varphi (e_j)^2 \| & \le & \| \varphi (e_j) - \varphi (\be_{M_r}) \varphi (e_j) \| + \| (\varphi (\be_{M_r}) - \varphi (e_j) )\varphi (e_j) \| \\
& < & \gamma^{\halb} + 0 \; ,
\end{eqnarray*}
so again the $\varphi (e_j)$ are close to (pairwise orthogonal) projections $p_j$.

Set $c_j := (p_j \varphi (e_j) p_j)^{-\halb} \in p_j A p_j$, then $c_j \varphi (e_j) c_j = p_j$, $j = 1 , \ldots , r$. Again $c_i \perp c_j$, $i \neq j$. Furthermore
\[
\sum^{r}_1 p_j = \sum^{r}_1 g_{\halb} (\varphi (e_j)) = g_{\halb} \left( \sum^{r}_1 \varphi (e_j) \right) = g_{\halb} (\varphi (\be_{M_r})) = p \; ,
\]
and setting $d := \sum^{r}_1 c_j$, one has $d \varphi (\be_{M_r}) d = p$.

But $d$ and $\varphi (\be_{M_r})$ commute, so $d^2 = (p \varphi (\be_{M_r}) p)^{-1}$ and $d = c$. It follows that
\[
\varphi' (e_j) = c \varphi (e_j) c = d \varphi (e_j) d = c_j \varphi (e_j) c_j = p_j
\]
is a projection, and that $\varphi' (e_i) \perp \varphi' (e_j)$ if $i \neq j$.

The $e_j$ were arbitrary, therefore $\ord \varphi' = 0$. Furthermore, $\varphi' (\be_{M_r}) = \sum^{r}_1 p_j$ is a projection, so $\varphi'$ is a $*$-homomorphism by b).
\end{proof}\\

\subsection[{\rm Completely positive rank zero}]{\sc Completely positive rank zero}

In this section we will be concerned with the case where the completely positive rank is zero. As it turns out, $\cpr A = 0$ if and only if $A$ is an $AF$ algebra, i.e.\ a direct limit of finite-dimensional $C^*$-algebras.

\bn{\label{local-characterization}}
Recall Bratelli's local characterization of $AF$ algebras, which is independent of a special sequence of subalgebras (cf. \cite{Br} or \cite{Da}, Theorem 3.3.4):

\begin{theorem} 
A separable $C^*$-algebra $A$ is $AF$ if and only if it satisfies the following condition: 
  \begin{itemize}
  \item [($*$)] for all $\varepsilon > 0$ and $a_1 , \ldots , a_n \in A$ there is a finite-dimensional $C^*$-subalgebra $F \subset A$ such that $\dist (a_i , F) < \varepsilon$ for $i = 1 , \ldots , n$. 
  \end{itemize}
\end{theorem}
\en

\bn
It is well-known that a commutative $C^*$-algebra $A$ is $AF$ precisely when $\dim \hat{A} = 0$. In Theorem \ref{AF} we show that this remains true in the noncommutative case if one replaces the covering dimension of the spectrum by the completely positive rank of the algebra.

That only $AF$ algebras are zero-dimensional seems to be a very special feature of the completely positive rank (in comparison, for example, to the real rank or stable rank, cf. Section 5). 

However, one would expect from any dimension theory generalizing covering dimension to noncommutative $C^*$-algebras that it is zero on $AF$ algebras: This is because, topologically, a full matrix algebra might be viewed as a point (or, more sophisticated, as $n$ equivalent points); so an $AF$ algebra is something like a limit of finite discrete (noncommutative) spaces. Thus if a dimension theory has sufficiently nice properties, it should be zero on $AF$ algebras.
\en\\

\bn{\label{AF}}
\begin{theorem}
For a nonzero separable $C^*$-algebra $A$ the following are equivalent:\\
a) $A$ is $AF$\\
b) $\cpr A = 0$.
\end{theorem}

\begin{proof}
  a) $\Rightarrow$ b): $\cpr F = 0$ for every finite-dimensional $C^*$-algebra $F$. By Proposition \ref{limits} $\cpr A = 0$ for an $AF$ algebra $A$.

b) $\Rightarrow$ a): The idea of the proof is simple: We have just seen that a unital c.p.\ map $\varphi : F \to A$ of strict order zero in fact is a $*$-homomorphism. If it is almost unital, it is close to a $*$-homomorphism. Problems only arise when $A$ has no unit. But in this case $F$ may be chosen to contain a sufficiently large hereditary subalgebra on which $\varphi$ is close to a $*$-homomorphism. Thus we can find enough finite-dimensional subalgebras of $A$; using the local characterization of $AF$ algebras (Theorem \ref{local-characterization}) this implies that $A$ in fact is $AF$.

Take $a_1 , \ldots , a_k \in A_+ , \, \delta > 0 , \, \varepsilon  > 0 , \, (u_{\lambda})_{\Lambda}$ a positive approximate unit for $A$ with $ \|a_i \| , \| u_{\lambda} \| \le 1$. By making $\varepsilon$ smaller if necessary, we may assume $\varepsilon^{1/4} < \gamma$, where we obtain $\gamma$ from $\delta$ by Proposition \ref{ord=0} c).

If we choose $u_{\lambda_0} , \,  u_{\lambda_1}$ with $\|u_{\lambda_0} u_{\lambda_1} - u_{\lambda_1} \|$ and $\, \| u_{\lambda_j} a_i - a_i \| ,\, j = 0,1$, sufficiently small, there is a c.p.\ approximation $(F , \psi , \varphi)$, with $\ord \varphi = 0$ and the following properties: 

\noindent $ \left.
\begin{array}{ll}
\rm (i) & \| \varphi \psi (a_i) - a_i \| < \varepsilon \; , \\
\rm (ii) & \| \varphi \psi (u_{\lambda_j}) - u_{\lambda_j} \| < \varepsilon \; , \\
\rm (iii) & \| \varphi \psi (u_{\lambda_j}) a_i - a_i \| < \varepsilon \; , 
\end{array} \right\} i = 1 , \ldots , k , \; j = 0,1$\\
$\begin{array}{ll}
\rm (iv) & \| u_ {\lambda_1} - h u_{\lambda_1} \| < \varepsilon \; , \; \mbox{where} \; h := \varphi (\be_F) \; , \\
\rm (v) & \| \varphi \psi (u_{\lambda_1}) - h \varphi \psi (u_{\lambda_1})\| < \varepsilon \; , \\
\rm (vi) & \varphi (p) \perp \varphi (q) \; \mbox{for orthogonal projections} \; p,q \in F\; .
\end{array}$

For (i), (ii), (iii) we approximate the $a_i$ and $u_{\lambda_j}$ sufficiently well by $(F , \psi , \varphi)$, for (iv) and (v) in addition we use Proposition \ref{almost-units}; (vi) follows from the fact, that $\ord \varphi = 0$. 

We have $\psi (u_{\lambda_1}) = \sum^r_1 \mu_k e_k$, where $0 \le \mu_k \le 1$ and $\{ e_1 , \ldots , e_r \} \subset F$ is elementary. Set $q := \sum_{\mu_k \ge \sqrt{\varepsilon}} e_k$, then
\[
\varphi (q) = \sum_{\mu_k \ge \sqrt{\varepsilon}} \varphi (e_k) \le \frac{1}{\sqrt{\varepsilon}} \varphi \psi (u_{\lambda_1}) \; ,
\]
and 
\[
\varphi (q) \ge \varphi \psi (u_{\lambda_1}) - \sqrt{\varepsilon} \cdot \be \; .
\]
Here, $\be$ denotes the unit of $A^+$. Then
\begin{eqnarray*}
  \| \varphi (q) - h \varphi (q) \|^2 & = & \| (\be-h) \varphi (q)^2 (\be -h)\| \le \| (\be-h) \varphi (q) (\be-h) \| \\
& \le & \frac{1}{\sqrt{\varepsilon}} \| (\be-h) \varphi \psi (u_{\lambda_1}) (\be -h) \| < \frac{\varepsilon}{\sqrt{\varepsilon}} = \sqrt{\varepsilon} \; ,
\end{eqnarray*}
so $\| \varphi (q) - h \varphi (q) \| < \varepsilon^{\frac{1}{4}}$.

We now have
\begin{eqnarray*}
  \| \varphi (q) - \varphi (q)^2 \| & \le & \| \varphi (q) - \varphi (q) h\| + \|\varphi (q) (h - \varphi (q)) \| \\
& = & \| \varphi (q) - \varphi (q) h \| + \| \varphi (q) \varphi (\be_F - q)\| \\
& \stackrel{\rm (vi)}{<} & \varepsilon^{\frac{1}{4}} + 0 \le \gamma \; .
\end{eqnarray*}

Now define $F' := q F q$; then $q = \be_{F'}$ and by Proposition \ref{ord=0} b) there is a $*$-homomorphism $\varphi' : F' \to A$ such that $\| \varphi' - \varphi|_{F'} \| < \delta$. Clearly, $ \varphi' (F') \subset A$ is a finite-dimensional $C^*$-subalgebra. 

It remains to be shown that $a_i$ is close to $\varphi' (F')$ for $i = 1 , \ldots , k$. We have
\begin{eqnarray*}
  \| a_i - a_i \varphi (q) \|^2 & = & \| a_i (\be - \varphi (q))^2 a_i \| \\
& \le & \| a_i (\be - \varphi (q)) a_i \| \\
& \le & \| a_i ((1 + \sqrt{\varepsilon}) \cdot \be - \varphi \psi (u_{\lambda_1})) a_i \| \\
& \stackrel{\rm (iii)}{<} & \sqrt{\varepsilon} + \varepsilon < 2 \sqrt{\varepsilon} \; .
\end{eqnarray*}
Using $\| \varphi (q) - \varphi (q)^2\| < \varepsilon^{\frac{1}{4}}$ and Proposition \ref{multiplicativity} we obtain 
\[
\| \varphi (q \psi (a_i) q) - \varphi (q) \varphi \psi (a_i) \varphi (q) \| < 2 \varepsilon^{\frac{1}{8}}.
\]
Therefore
\begin{eqnarray*}
  \| a_i - \varphi (q \psi (a_i) q)\| & \le & \| a_i - \varphi (q) a_i \varphi (q) \| + \| \varphi (q) a_i \varphi (q) - \varphi (q) \varphi \psi (a_i) \varphi (q) \| \\
& & + \| \varphi (q) \varphi \psi (a_i) \varphi (q) - \varphi (q \psi (a_i) q)\| \\
& < & 2 \cdot \sqrt{2} \varepsilon^{\frac{1}{4}} + \varepsilon + 2 \varepsilon^{\frac{1}{8}} \; .
\end{eqnarray*}
Furthermore $\| a_i - \varphi' (q \psi (a_i) q) \|$ is small, because $\| \varphi' - \varphi \| < \delta$, so $a_i$ is close to $\varphi' (F')$.

But then by Theorem \ref{local-characterization}, $A$ is an $AF$ algebra.
\end{proof}
\en

\section[{\sc Other concepts}]{\large \sc Other concepts of noncommutative covering dimension}

There are several other concepts of noncommutative covering dimension, among which the stable and the real rank (introduced by Rieffel and by Brown and Pedersen, respectively) are of special interest. Below we describe these notions and some of their properties; furthermore we briefly comment on Lin's tracial rank and on Murphy's analytic rank.
All these (definitely different) concepts of covering dimension in some sense illustrate various points of view on noncommutative topology, which are very closely related in the commutative case.

\subsection[{\rm Stable and real rank}]{\sc Stable and real rank}

\bn
\begin{defn}{\rm (cf.\ \cite{Ri1}, Definition 1.4 and \cite{BP}, Section 1)}
Let $A$ be a unital $C^*$-algebra.\\
(i) The stable rank of $A$ is less than or equal to $n$, $\sr A \le n$, if the following holds:
\begin{itemize}
\item [($*$)] For each $(n+1)$-tuple $(x_0 , \ldots , x_n) \in A^{n+1}$ and every $\varepsilon > 0$ there is an $(n+1)$-tuple $(y_0 , \ldots , y_n) \in A^{n+1}$ such that $\sum y^*_k y_k$ is invertible and $\sum \| x_k - y_k \| < \varepsilon$.
\end{itemize}
(ii) The real rank of $A$ is less than or equal to $n$, $\rr A \le n$, if ($*$) holds with the additional condition that the $x_i$ and $y_i$ are in $A_{sa}$.\\
If $A$ is non-unital, define $\sr A := \sr A^+$ and $\rr A := \rr A^+$.   
\end{defn}
\en

\bn{\label{dim=rr}}
Whereas the notion of the completely positive rank is modelled after Definition \ref{covering-dimension}, the above definition is inspired by the characterization of covering dimension in \cite{Pe}, Theorem 3.3.2. Both concepts generalize covering dimension (cf.\ \cite{Ri1}, Proposition 1.7 and \cite{BP}, Proposition 1.1): If $X$ is a compact space, then $\sr (\Ch (X)) = [(\dim X)/2] + 1$ (where $[s]$ denotes the integer part of $s$) and $\rr (\Ch (X)) = \dim X$.

The stable rank does not directly generalize covering dimension. However, one of Rieffel's original motivations for introducing this notion was to obtain stability results for $C^*$-algebra $K$-theory analogous to those of topological $K$-theory. These classical results often involve terms of the form $[(\dim X)/ 2]+1$. So, from this point of view, the stable rank is an appropriate analogue of covering dimension.
\en

\bn
Stable rank behaves nicely with respect to ideals, quotients and inductive limits, whereas the behavior under taking matrix algebras is rather surprising (cf.\ \cite{Ri1}, Theorems 5.2 and 5.4): For any $r \in \N$ we have
\[
\sr (M_r (A)) = \{ (\sr A - 1) / r \} + 1 \; ,
\]
where $\{ s \}$ denotes the least integer greater than $s$. Furthermore we have $\sr (A \otimes \Kh) = 1$ if $\sr A = 1$; otherwise $\sr (A \otimes \Kh) = 2$. Also, it turns out that $\rr A \le 2 \cdot \sr A - 1$.
\en

\bn
\begin{examples} (cf.\ \cite{Ri1} and \cite{BP}) All $AF$ algebras have stable rank one and real rank zero; von Neumann algebras have real rank zero. The irrational rotation algebras and the Bunce--Deddens algebras have real rank zero and so do the Cuntz algebras $\Oh_n$ and, more generally, all simple purely infinite $C^*$-algebras (cf.\ \cite{Cu}). In contrast, every unital $C^*$-algebra containing two isometries with orthogonal range projections (in particular each simple unital purely infinite $C^*$-algebra) has infinite stable rank. In view of this it might be a little surprising that the Toeplitz algebra has stable rank two. If $\alpha$ is an automorphism of the $C^*$-algebra $A$, then $\sr (A \rtimes_{\alpha} \Z) \le \sr A + 1$.
\end{examples}
\en

\bn
It is natural to ask how the (strong) completely positive rank is related to the stable and real rank.

So far, in all our examples the real rank is less than or equal to the completely positive rank.

We have also seen that stable and real rank both break down under taking matrix algebras. The completely positive rank behaves very different in this respect, since, as we shall see in \cite{Wi}, $\cpr (M_r \otimes \Ch (X))= \cpr (\Ch (X)) (= \dim X)$. These observations suggest the following

\begin{conj}
  For any $C^*$-algebra $A$, $\rr A \le \cpr A$.
\end{conj}
\en

\subsection[{\rm Tracial rank}]{\sc Tracial rank}

For recent development in classification of nuclear $C^*$-algebras tracially $AF$ and tracially $AI$ algebras (introduced by Lin in \cite{Li1}) have been of certain interest. In \cite{Li3}, Lin has given a notion of covering dimension closely related to these concepts.

\bn
Denote by $I^{(k)}$ the class of all unital $C^*$-algebras which are hereditary subalgebras of $C^*$-algebras of the form $\Ch (X) \otimes M_r$, where $X$ is a $k$-dimensional finite $CW$-complex.\\
A unital $C^*$-algebra $A$ is said to have tracial rank less than or equal to $k$, $\tr A \le k$, if for any $\varepsilon > 0$ and any finite subset $G \subset A$ there is a $C^*$-subalgebra $B \in I^{(k)}$ of $A$ such that\\
(1) $\| [x , \be_B] \| < \varepsilon \; \forall x \in G$,\\
(2) $ \dist (\be_B x \be_B , B) < \varepsilon \; \forall x \in G$, \\
(3) $\be_A - \be_B$ is small in an appropriate sense.\\
The exact formulation of (3) uses comparison theory and is a little technical.
\en

\bn
For a compact space $X$, $\tr (\Ch (X)) = \dim X$. Furthermore, if $J \lhd A$ is an ideal, then $\tr (A / J) \le \tr A$; if $B \subset_{\her} A$ is a hereditary $C^*$-subalgebra with unit, then $\tr B \le \tr A$.
\en

\bn
It turns out that irrational rotation algebras and Bunce--Deddens algebras have tracial rank zero. On the other hand, Blackadar's unital projectionless $C^*$-algebra has infinite tracial rank. These examples in particular show that neither the tracial rank nor the completely positive rank dominates the other and at the present stage we do not know how these concepts might be related to each other.
\en

\subsection[{\rm Analytic rank}]{\sc Analytic rank}

In \cite{Mu2} Murphy introduced the analytic rank, another notion of covering dimension for $C^*$-algebras, the topological analogue of which is defined directly in terms of the continuous functions on a space. Below we give the definition and some sample results (see \cite{Mu2} for details).

\bn
\begin{defn}
(i) Let $A$ be a unital $C^*$-algebra, $B \subset A$ a $C^*$-subalgebra. $B$ is said to be analytic, if $\be_A \in B$ and if for each $a \in A_{sa}$ with $a^2 \in B$ automatically $a \in B$.\\
For a subset $S \subset A$ denote by $A[S]$ the smallest analytic $C^*$-subalgebra containing $S$. If $S \subset A_{sa}$ and $A[S] = A$, then $S$ is said to be an analytic base of $A$.\\
(ii) The analytic rank of $A$ does not exceed $n$, $\rk A \le n$, if $A$ has an analytic base consisting of $n$ elements.
\end{defn}
\en

\bn
Note that, in the commutative case, this definition coincides with the ``classical'' analytical rank, cf.\ \cite{Pe}, Ch.\ 10.4. If $X$ is a compact metrizable space, then $\rk (\Ch (X)) = \dim X$ (\cite{Pe}, Proposition 10.4.21).
\en

\bn
If $(A_{\lambda})_{\Lambda}$ is a family of unital subalgebras of $A$ such that $\bigcup_{\Lambda} A_{\lambda}$ generates $A$, then $\rk A \le \sum_{\Lambda} \rk (A_{\lambda})$. As a consequence, $\rk (A_1 \otimes A_2) \le \rk (A_1) + \rk (A_2)$. Similar results for the other notions of noncommutative covering dimension, it seems, would be hard to obtain.

If $(A,G,\alpha)$ is a $C^*$-dynamical system with $G$ a countable discrete abelian group, then $\rk (A \rtimes_{\alpha} G) \le \rk A + \dim \hat{G}$.

For any unital $C^*$-algebra and any $r > 1$, $\rk (M_r(A)) =0$; in particular, analytic rank zero does not pass to hereditary subalgebras, whereas it is preserved under inductive limits.

If $\rr A =0$, then $\rk A =0$. If $A$ is generated by $n$ (partial) isometries, then $\rk A \le n$, so for the Toeplitz algebra $\Th$ we have $\rk \Th = 1$ (the rank cannot be zero since $\Th$ has a one-dimensional quotient). For the full group $C^*$-algebra $C^*(\F_n)$ we obtain $\rk (C^*(\F_n)) = n$.

Note that, despite its rather strange behavior in some respects, e.g.\ under taking matrix algebras, the analytic rank can have reasonable values also for non-nuclear $C^*$-algebras, like $C^*(\F_n)$.

Although for all our examples the analytic rank is smaller than the completely positive rank, we do not at all have a general result.
\en

\section[{\sc Some open questions}]{\large \sc Some open questions}

\altbn
For locally compact spaces $X$ and $Y$, we have
\[
\cpr (\Ch_0(X \times Y)) \le \cpr(\Ch_0(X)) + \cpr(\Ch_0(Y)).
\]
One might expect a similar behavior for more general $C^*$-algebras. However, it seems a little daring to write down a formula like
\[
\cpr (A \otimes B) \le \cpr A + \cpr B,
\]
for at the present stage the almost only body of evidence is the lack of counterexamples. We do not even have a satisfactory answer in the cases where one of the factors is commutative or just a matrix algebra, although it seems likely to be true that $\cpr (M_r(A)) \le \cpr A$; this indeed is the case if $\cpr A \le 1$ or if $A$ is simple, as we will show in \cite{Wi}. But then still the question remains in which cases we even have equality.

\alten

\altbn
Using the countable sum theorem \ref{countable-sums} it is not hard to see that, for an open subset $U$ of a second countable locally compact space $X$, we have $\dim U \le \dim X$. This provokes the following question: If $A$ is a $C^*$-algebra and $B \subset_{\her} A$, do we have $\cpr B \le \cpr A$? In \cite{Wi} we show that this is true at least if $B$ is an ideal of $A$. 

\alten

\altbn
So far all examples for which we were able to determine the completely positive rank, were (stably) finite. It is natural to ask wether infinite $C^*$-algebras can have finite completely positive rank. In particular, what values do we obtain for the Toeplitz algebra and for the Cuntz algebras? \\
One can describe concrete systems of c.p.\ approximations for the $\Oh_n$ and for $\Th$, however, we did not succeed in finding one where the order of the $\varphi$'s is bounded.
 
\alten

\altbn
If $G$ is an amenable group, then the natural surjection $C^*(G) \to C^*_r(G)$ is an isomorphism and the group $C^*$-algebra is nuclear. Unfortunately, it is not clear how to obtain c.p.\ approximations for $C^*(G)$ systematically (depending, for example, on a left invariant mean on $G$), so we do not have general results on the completely positive rank of group $C^*$-algebras. Of course, if $G$ is abelian, then $C^*(G) = \Ch_0(\hat{G})$, so $\cpr (C^*(G)) = \dim \hat{G}$.

\alten

\altbn
If $(A,G,\alpha)$ is a $C^*$-dynamical system with $A$ nuclear and $G$ amenable, then the crossed product $A \rtimes_{\alpha} G$ is nuclear, but again we do not know how to obtain c.p.\ approximations in a systematic way.\\
If $G$ is abelian, then $A \rtimes_{\alpha} G$ is often regarded as a generalized ``skew'' tensor product of $A$ and $C^*(G) = \Ch_0 (\hat{G})$. One might ask if there is a formula like 
\[
\cpr (A \rtimes_{\alpha} G) \le \cpr A + \dim \hat{G}.
\]
This is not true, for Blackadar in \cite{Bl3} has given an example of a crossed product of an $AF$ algebra with a $\Z_2$-action, which is not an $AF$ algebra. Note that, because of a different counterexample, the above estimate also fails for the real rank (cf.\ \cite{BP}).\\
In the case of the irrational rotation algebras we know the completely positive rank because $A_{\theta}$ can be written as a limit of circle algebras. More generally, consider $\Ch(M) \rtimes_{\alpha} \Z$ for some smooth manifold $M$ and a minimal diffeomorphism $\alpha$ of $M$. Using a characterization of crossed products of this form recently given by Q.\ Lin and N.\ C.\ Phillips it seems to be possible to show that 
\[
\cpr(\Ch(M) \rtimes_{\alpha} \Z) \le \dim M \, .
\] 
\alten

\altbn
Covering dimension is a topological invariant which is defined via intersections of the members of open coverings, and in this respect it is related to \v{C}ech (co)homology (cf.\ \cite{ES}). One of the original reasons for introducing yet another notion of noncommutative covering dimension was it, to find out how something like a \v{C}ech cohomology for nuclear $C^*$-algebras (based on c.p.\ approximations) could possibly be constructed.\\
Of course one should not expect such a theory to be too well-behaved. In \cite{Cu2}, Cuntz has shown that a covariant functor from $C^*$-algebras into abelian groups with sufficiently nice abstract properties cannot avoid being $K$-theory (at least on a large subclass of $C^*$-algebras). The main point of the argument is, that such a functor must be Bott periodic, which in turn is based on a sophisticated use of the Toeplitz algebra.\\
There still might be some hope to obtain a notion of \v{C}ech cohomology at least for (a subclass of) finite $C^*$-algebras. But even then it follows from Loring's $AF$-embeddings of the rational rotation algebras (cf.\ \cite{Lo}), that such a functor, if it has some (rather weak) abstract properties, cannot be too well-behaved either. In particular, it will not admit a Chern character which restricts to filtration by dimension in the commutative case.

\alten


\tableofcontents

\clearpage
\thispagestyle{empty}
\vspace*{\fill}

\clearpage
\thispagestyle{empty}
\vspace*{\fill}


\end{document}